\magnification=\magstephalf
\magnification=\magstephalf
\def\temp{1.35}%
\let\tempp=\relax
\expandafter\ifx\csname psboxversion\endcsname\relax
  \message{PSBOX(\temp)}%
\else
    \ifdim\temp cm>\psboxversion cm
      \message{PSBOX(\temp)}%
    \else
      \message{PSBOX(\psboxversion) is already loaded: I won't load
        PSBOX(\temp)!}%
      \let\temp=\psboxversion
      \let\tempp= 
    \fi
\fi
\tempp
\message{by Jean Orloff: loading ...}
\let\psboxversion=\temp
\catcode`\@=11
%
%
\def\psfortextures{
\def\PSspeci@l##1##2{%
\special{illustration ##1\space scaled ##2}%
}}%
\def\psfordvitops{
\def\PSspeci@l##1##2{%
\special{dvitops: import ##1\space \the\drawingwd \the\drawinght}%
}}%
\def\psfordvips{
\def\PSspeci@l##1##2{%
\d@my=0.1bp \d@mx=\drawingwd \divide\d@mx by\d@my
\includegraphics{##1\space}}}%
\def\psforoztex{
\def\PSspeci@l##1##2{%
\special{##1 \space
      ##2 1000 div dup scale
      \number-\psllx\space\space \number-\pslly\space\space translate
}}}%
\def\psfordvitps{
\def\dvitpsLiter@ldim##1{\dimen0=##1\relax
\special{dvitps: Literal "\number\dimen0\space"}}%
\def\PSspeci@l##1##2{%
\at(0bp;\drawinght){%
\special{dvitps: Include0 "psfig.psr"}
\dvitpsLiter@ldim{\drawingwd}%
\dvitpsLiter@ldim{\drawinght}%
\dvitpsLiter@ldim{\psllx bp}%
\dvitpsLiter@ldim{\pslly bp}%
\dvitpsLiter@ldim{\psurx bp}%
\dvitpsLiter@ldim{\psury bp}%
\special{dvitps: Literal "startTexFig"}%
\special{dvitps: Include1 "##1"}%
\special{dvitps: Literal "endTexFig"}%
}}}%
\def\psfordvialw{
\def\PSspeci@l##1##2{
\special{language "PostScript",
position = "bottom left",
literal "  \psllx\space \pslly\space translate
  ##2 1000 div dup scale
  -\psllx\space -\pslly\space translate",
include "##1"}
}}%
\def\psforptips{
\def\PSspeci@l##1##2{{
\d@mx=\psurx bp
\advance \d@mx by -\psllx bp
\divide \d@mx by 1000\multiply\d@mx by \xscale
\incm{\d@mx}
\let\tmpx\dimincm
\d@my=\psury bp
\advance \d@my by -\pslly bp
\divide \d@my by 1000\multiply\d@my by \xscale
\incm{\d@my}
\let\tmpy\dimincm
\d@mx=-\psllx bp
\divide \d@mx by 1000\multiply\d@mx by \xscale
\d@my=-\pslly bp
\divide \d@my by 1000\multiply\d@my by \xscale
\at(\d@mx;\d@my){\special{ps:##1 x=\tmpx cm, y=\tmpy cm}}
}}}%
\def\psonlyboxes{
\def\PSspeci@l##1##2{%
\at(0cm;0cm){\boxit{\vbox to\drawinght
  {\vss\hbox to\drawingwd{\at(0cm;0cm){\hbox{({\tt##1})}}\hss}}}}
}}%
\def\psloc@lerr#1{%
\let\savedPSspeci@l=\PSspeci@l%
\def\PSspeci@l##1##2{%
\at(0cm;0cm){\boxit{\vbox to\drawinght
  {\vss\hbox to\drawingwd{\at(0cm;0cm){\hbox{({\tt##1}) #1}}\hss}}}}
\let\PSspeci@l=\savedPSspeci@l
}}%
%
%
\newread\pst@mpin
\newdimen\drawinght\newdimen\drawingwd
\newdimen\psxoffset\newdimen\psyoffset
\newbox\drawingBox
\newcount\xscale \newcount\yscale \newdimen\pscm\pscm=1cm
\newdimen\d@mx \newdimen\d@my
\newdimen\pswdincr \newdimen\pshtincr
\let\ps@nnotation=\relax
{\catcode`\|=0 |catcode`|\=12 |catcode`|
|catcode`#=12 |catcode`*=14
|xdef|backslashother{\}*
|xdef|percentother{
|xdef|tildeother{~}*
|xdef|sharpother{#}*
}%
\def\R@moveMeaningHeader#1:->{}%
\def\uncatcode#1{%
\edef#1{\expandafter\R@moveMeaningHeader\meaning#1}}%
\def\execute#1{#1}
\def\psm@keother#1{\catcode`#112\relax}
\def\executeinspecs#1{%
\execute{\begingroup\let\do\psm@keother\dospecials\catcode`\^^M=9#1\endgroup}}%
\def\@mpty{}%
\def\matchexpin#1#2{
  \fi%
  \edef\tmpb{{#2}}%
  \expandafter\makem@tchtmp\tmpb%
  \edef\tmpa{#1}\edef\tmpb{#2}%
  \expandafter\expandafter\expandafter\m@tchtmp\expandafter\tmpa\tmpb\endm@tch%
  \if\match%
}%
\def\matchin#1#2{%
  \fi%
  \makem@tchtmp{#2}%
  \m@tchtmp#1#2\endm@tch%
  \if\match%
}%
\def\makem@tchtmp#1{\def\m@tchtmp##1#1##2\endm@tch{%
  \def\tmpa{##1}\def\tmpb{##2}\let\m@tchtmp=\relax%
  \ifx\tmpb\@mpty\def\match{YN}%
  \else\def\match{YY}\fi%
}}%
\def\incm#1{{\psxoffset=1cm\d@my=#1
 \d@mx=\d@my
  \divide\d@mx by \psxoffset
  \xdef\dimincm{\number\d@mx.}
  \advance\d@my by -\number\d@mx cm
  \multiply\d@my by 100
 \d@mx=\d@my
  \divide\d@mx by \psxoffset
  \edef\dimincm{\dimincm\number\d@mx}
  \advance\d@my by -\number\d@mx cm
  \multiply\d@my by 100
 \d@mx=\d@my
  \divide\d@mx by \psxoffset
  \xdef\dimincm{\dimincm\number\d@mx}
}}%
%
\newif\ifNotB@undingBox
\newhelp\PShelp{Proceed: you'll have a 5cm square blank box instead of
your graphics.}%
\def\s@tsize#1 #2 #3 #4\@ndsize{
  \def\psllx{#1}\def\pslly{#2}%
  \def\psurx{#3}\def\psury{#4}
  \ifx\psurx\@mpty\NotB@undingBoxtrue
  \else
    \drawinght=#4bp\advance\drawinght by-#2bp
    \drawingwd=#3bp\advance\drawingwd by-#1bp
  \fi
  }%
\def\sc@nBBline#1:#2\@ndBBline{\edef\p@rameter{#1}\edef\v@lue{#2}}%
\def\g@bblefirstblank#1#2:{\ifx#1 \else#1\fi#2}%
{\catcode`\%=12
\xdef\B@undingBox{
\def\ReadPSize#1{
 \readfilename#1\relax
 \let\PSfilename=\lastreadfilename
 \openin\pst@mpin=#1\relax
 \ifeof\pst@mpin \errhelp=\PShelp
   \errmessage{I haven't found your postscript file (\PSfilename)}%
   \psloc@lerr{was not found}%
   \s@tsize 0 0 142 142\@ndsize
   \closein\pst@mpin
 \else
   \if\matchexpin{\GlobalInputList}{, \lastreadfilename}%
   \else\xdef\GlobalInputList{\GlobalInputList, \lastreadfilename}%
     \immediate\write\psbj@inaux{\lastreadfilename,}%
   \fi%
   \loop
     \executeinspecs{\catcode`\ =10\global\read\pst@mpin to\n@xtline}%
     \ifeof\pst@mpin
       \errhelp=\PShelp
       \errmessage{(\PSfilename) is not an Encapsulated PostScript File:
           I could not find any \B@undingBox: line.}%
       \edef\v@lue{0 0 142 142:}%
       \psloc@lerr{is not an EPSFile}%
       \NotB@undingBoxfalse
     \else
       \expandafter\sc@nBBline\n@xtline:\@ndBBline
       \ifx\p@rameter\B@undingBox\NotB@undingBoxfalse
         \edef\t@mp{%
           \expandafter\g@bblefirstblank\v@lue\space\space\space}%
         \expandafter\s@tsize\t@mp\@ndsize
       \else\NotB@undingBoxtrue
       \fi
     \fi
   \ifNotB@undingBox\repeat
   \closein\pst@mpin
 \fi
\message{#1}%
}%
%
%
\def\psboxto(#1;#2)#3{\vbox{%
   \ReadPSize{#3}%
   \advance\pswdincr by \drawingwd
   \advance\pshtincr by \drawinght
   \divide\pswdincr by 1000
   \divide\pshtincr by 1000
   \d@mx=#1
   \ifdim\d@mx=0pt\xscale=1000
         \else \xscale=\d@mx \divide \xscale by \pswdincr\fi
   \d@my=#2
   \ifdim\d@my=0pt\yscale=1000
         \else \yscale=\d@my \divide \yscale by \pshtincr\fi
   \ifnum\yscale=1000
         \else\ifnum\xscale=1000\xscale=\yscale
                    \else\ifnum\yscale<\xscale\xscale=\yscale\fi
              \fi
   \fi
   \divide\drawingwd by1000 \multiply\drawingwd by\xscale
   \divide\drawinght by1000 \multiply\drawinght by\xscale
   \divide\psxoffset by1000 \multiply\psxoffset by\xscale
   \divide\psyoffset by1000 \multiply\psyoffset by\xscale
   \global\divide\pscm by 1000
   \global\multiply\pscm by\xscale
   \multiply\pswdincr by\xscale \multiply\pshtincr by\xscale
   \ifdim\d@mx=0pt\d@mx=\pswdincr\fi
   \ifdim\d@my=0pt\d@my=\pshtincr\fi
   \message{scaled \the\xscale}%
 \hbox to\d@mx{\hss\vbox to\d@my{\vss
   \global\setbox\drawingBox=\hbox to 0pt{\kern\psxoffset\vbox to 0pt{%
      \kern-\psyoffset
      \PSspeci@l{\PSfilename}{\the\xscale}%
      \vss}\hss\ps@nnotation}%
   \global\wd\drawingBox=\the\pswdincr
   \global\ht\drawingBox=\the\pshtincr
   \global\drawingwd=\pswdincr
   \global\drawinght=\pshtincr
   \baselineskip=0pt
   \copy\drawingBox
 \vss}\hss}%
  \global\psxoffset=0pt
  \global\psyoffset=0pt
  \global\pswdincr=0pt
  \global\pshtincr=0pt 
  \global\pscm=1cm 
}}%
%
%
\def\psboxscaled#1#2{\vbox{%
  \ReadPSize{#2}%
  \xscale=#1
  \message{scaled \the\xscale}%
  \divide\pswdincr by 1000 \multiply\pswdincr by \xscale
  \divide\pshtincr by 1000 \multiply\pshtincr by \xscale
  \divide\psxoffset by1000 \multiply\psxoffset by\xscale
  \divide\psyoffset by1000 \multiply\psyoffset by\xscale
  \divide\drawingwd by1000 \multiply\drawingwd by\xscale
  \divide\drawinght by1000 \multiply\drawinght by\xscale
  \global\divide\pscm by 1000
  \global\multiply\pscm by\xscale
  \global\setbox\drawingBox=\hbox to 0pt{\kern\psxoffset\vbox to 0pt{%
     \kern-\psyoffset
     \PSspeci@l{\PSfilename}{\the\xscale}%
     \vss}\hss\ps@nnotation}%
  \advance\pswdincr by \drawingwd
  \advance\pshtincr by \drawinght
  \global\wd\drawingBox=\the\pswdincr
  \global\ht\drawingBox=\the\pshtincr
  \global\drawingwd=\pswdincr
  \global\drawinght=\pshtincr
  \baselineskip=0pt
  \copy\drawingBox
  \global\psxoffset=0pt
  \global\psyoffset=0pt
  \global\pswdincr=0pt
  \global\pshtincr=0pt 
  \global\pscm=1cm
}}%
%
\def\psbox#1{\psboxscaled{1000}{#1}}%
\newif\ifn@teof\n@teoftrue
\newif\ifc@ntrolline
\newif\ifmatch
\newread\j@insplitin
\newwrite\j@insplitout
\newwrite\psbj@inaux
\immediate\openout\psbj@inaux=psbjoin.aux
\immediate\write\psbj@inaux{\string\joinfiles}%
\immediate\write\psbj@inaux{\jobname,}%
%
%
\def\toother#1{\ifcat\relax#1\else\expandafter%
  \toother@ux\meaning#1\endtoother@ux\fi}%
\def\toother@ux#1 #2#3\endtoother@ux{\def\tmp{#3}%
  \ifx\tmp\@mpty\def\tmp{#2}\let\next=\relax%
  \else\def\next{\toother@ux#2#3\endtoother@ux}\fi%
\next}%
%
%
\let\readfilenamehook=\relax
\def\re@d{\expandafter\re@daux}
\def\re@daux{\futurelet\nextchar\stopre@dtest}%
\def\re@dnext{\xdef\lastreadfilename{\lastreadfilename\nextchar}%
  \afterassignment\re@d\let\nextchar}%
\def\stopre@d{\egroup\readfilenamehook}%
\def\stopre@dtest{%
  \ifcat\nextchar\relax\let\nextread\stopre@d
  \else
    \ifcat\nextchar\space\def\nextread{%
      \afterassignment\stopre@d\chardef\nextchar=`}%
    \else\let\nextread=\re@dnext
      \toother\nextchar
      \edef\nextchar{\tmp}%
    \fi
  \fi\nextread}%
\def\readfilename{\bgroup%
  \let\\=\backslashother \let\%=\percentother \let\~=\tildeother
  \let\#=\sharpother \xdef\lastreadfilename{}%
  \re@d}%
%
%
\xdef\GlobalInputList{\jobname}%
\def\psnewinput{%
  \def\readfilenamehook{
    \if\matchexpin{\GlobalInputList}{, \lastreadfilename}%
    \else\xdef\GlobalInputList{\GlobalInputList, \lastreadfilename}%
      \immediate\write\psbj@inaux{\lastreadfilename,}%
    \fi%
    \let\readfilenamehook=\relax%
    \ps@ldinput\lastreadfilename\relax%
  }\readfilename%
}%
\expandafter\ifx\csname @@input\endcsname\relax    
  \immediate\let\ps@ldinput=\input\def\input{\psnewinput}%
\else
  \immediate\let\ps@ldinput=\@@input
  \def\@@input{\psnewinput}%
\fi%
\def\nowarnopenout{%
 \def\warnopenout##1##2{%
   \readfilename##2\relax
   \message{\lastreadfilename}%
   \immediate\openout##1=\lastreadfilename\relax}}%
\def\warnopenout#1#2{%
 \readfilename#2\relax
 \def\t@mp{TrashMe,psbjoin.aux,psbjoint.tex,}\uncatcode\t@mp
 \if\matchexpin{\t@mp}{\lastreadfilename,}%
 \else
   \immediate\openin\pst@mpin=\lastreadfilename\relax
   \ifeof\pst@mpin
     \else
     \edef\tmp{{If the content of this file is precious to you, this
is your last chance to abort (ie press x or e) and rename it before
retexing (\jobname). If you're sure there's no file
(\lastreadfilename) in the directory of (\jobname), then go on: I'm
simply worried because you have another (\lastreadfilename) in some
directory I'm looking in for inputs...}}%
     \errhelp=\tmp
     \errmessage{I may be about to replace your file named \lastreadfilename}%
   \fi
   \immediate\closein\pst@mpin
 \fi
 \message{\lastreadfilename}%
 \immediate\openout#1=\lastreadfilename\relax}%
{\catcode`\%=12\catcode`\*=14
\gdef\splitfile#1{*
 \readfilename#1\relax
 \immediate\openin\j@insplitin=\lastreadfilename\relax
 \ifeof\j@insplitin
   \message{! I couldn't find and split \lastreadfilename!}*
 \else
   \immediate\openout\j@insplitout=TrashMe
   \message{< Splitting \lastreadfilename\space into}*
   \loop
     \ifeof\j@insplitin
       \immediate\closein\j@insplitin\n@teoffalse
     \else
       \n@teoftrue
       \executeinspecs{\global\read\j@insplitin to\spl@tinline\expandafter
         \ch@ckbeginnewfile\spl@tinline
       \ifc@ntrolline
       \else
         \toks0=\expandafter{\spl@tinline}*
         \immediate\write\j@insplitout{\the\toks0}*
       \fi
     \fi
   \ifn@teof\repeat
   \immediate\closeout\j@insplitout
 \fi\message{>}*
}*
\gdef\ch@ckbeginnewfile#1
 \def\t@mp{#1}*
 \ifx\@mpty\t@mp
   \def\t@mp{#3}*
   \ifx\@mpty\t@mp
     \global\c@ntrollinefalse
   \else
     \immediate\closeout\j@insplitout
     \warnopenout\j@insplitout{#2}*
     \global\c@ntrollinetrue
   \fi
 \else
   \global\c@ntrollinefalse
 \fi}*
\gdef\joinfiles#1\into#2{*
 \message{< Joining following files into}*
 \warnopenout\j@insplitout{#2}*
 \message{:}*
 {*
 \edef\w@##1{\immediate\write\j@insplitout{##1}}*
\w@{
\w@{
\w@{
\w@{
\w@{
\w@{
\w@{
\w@{
\w@{
\w@{
\w@{\string\input\space psbox.tex}*
\w@{\string\splitfile{\string\jobname}}*
\w@{\string\let\string\autojoin=\string\relax}*
}*
 \expandafter\tre@tfilelist#1, \endtre@t
 \immediate\closeout\j@insplitout
 \message{>}*
}*
\gdef\tre@tfilelist#1, #2\endtre@t{*
 \readfilename#1\relax
 \ifx\@mpty\lastreadfilename
 \else
   \immediate\openin\j@insplitin=\lastreadfilename\relax
   \ifeof\j@insplitin
     \errmessage{I couldn't find file \lastreadfilename}*
   \else
     \message{\lastreadfilename}*
     \immediate\write\j@insplitout{
     \executeinspecs{\global\read\j@insplitin to\oldj@ininline}*
     \loop
       \ifeof\j@insplitin\immediate\closein\j@insplitin\n@teoffalse
       \else\n@teoftrue
         \executeinspecs{\global\read\j@insplitin to\j@ininline}*
         \toks0=\expandafter{\oldj@ininline}*
         \let\oldj@ininline=\j@ininline
         \immediate\write\j@insplitout{\the\toks0}*
       \fi
     \ifn@teof
     \repeat
   \immediate\closein\j@insplitin
   \fi
   \tre@tfilelist#2, \endtre@t
 \fi}*
}%
\def\autojoin{%
 \immediate\write\psbj@inaux{\string\into{psbjoint.tex}}%
 \immediate\closeout\psbj@inaux
 \expandafter\joinfiles\GlobalInputList\into{psbjoint.tex}%
}%
%
%
%
\def\centinsert#1{\ifhmode\vadjust\fi{%
  \midinsert\smallskip\line{\hss#1\hss}\smallskip\endinsert}}%
\def\psannotate#1#2{\vbox{%
  \def\ps@nnotation{#2\global\let\ps@nnotation=\relax}#1}}%
\newdimen\captwd \captwd=1.0\hsize
\def\pscaption#1#2{\vbox{%
   \setbox\drawingBox=#1
   \centerline{\copy\drawingBox}%
   \vskip\baselineskip
   \centerline{\vbox{\hsize=\captwd\setbox0=\hbox{#2}%
     \ifdim\wd0>\hsize
       \noindent\unhbox0\tolerance=5000
    \else\centerline{\box0}%
    \fi}%
}}}%
\def\psfig#1#2#3{\pscaption{\psannotate{#1}{#2}}{#3}}
\def\psfigurebox#1#2#3{\pscaption{\psannotate{\psbox{#1}}{#2}}{#3}}
%
\def\at(#1;#2)#3{\setbox0=\hbox{#3}\ht0=0pt\dp0=0pt
  \rlap{\kern#1\vbox to0pt{\kern-#2\box0\vss}}}%
%
\newdimen\gridht \newdimen\gridwd
\def\gridfill(#1;#2){%
  \setbox0=\hbox to 1\pscm
  {\vrule height1\pscm width.4pt\leaders\hrule\hfill}%
  \gridht=#1
  \divide\gridht by \ht0
  \multiply\gridht by \ht0
  \gridwd=#2
  \divide\gridwd by \wd0
  \multiply\gridwd by \wd0
  \advance \gridwd by \wd0
  \vbox to \gridht{\leaders\hbox to\gridwd{\leaders\box0\hfill}\vfill}}%
%
\def\fillinggrid{\at(0cm;0cm){\vbox{%
  \gridfill(\drawinght;\drawingwd)}}}%
%
%
\def\textleftof#1:{%
  \setbox1=#1
  \setbox0=\vbox\bgroup
    \advance\hsize by -\wd1 \advance\hsize by -2em}%
\def\textrightof#1:{%
  \setbox0=#1
  \setbox1=\vbox\bgroup
    \advance\hsize by -\wd0 \advance\hsize by -2em}%
\def\endtext{%
  \egroup
  \hbox to \hsize{\valign{\vfil##\vfil\cr%
\box0\cr%
\noalign{\hss}\box1\cr}}}%
%
\def\frameit#1#2#3{\hbox{\vrule width#1\vbox{%
  \hrule height#1\vskip#2\hbox{\hskip#2\vbox{#3}\hskip#2}%
        \vskip#2\hrule height#1}\vrule width#1}}%
\def\boxit#1{\frameit{0.4pt}{0pt}{#1}}%
\catcode`\@=12 
%
\psfortextures
 
\psfordvips    

\hsize = 15.5truecm
\vsize = 22truecm 




\nopagenumbers



\headline={\ifodd\pageno
           \rightheadline
            \else\leftheadline
           \fi
         }
\def\rightheadline{\pfont\hfil\number\pageno}
\def\leftheadline {\pfont\number\pageno\hfil}

\baselineskip=13pt
\parskip 2pt plus 1pt
\font\notefont=cmr9

\vbadness=10000
\widowpenalty=10000 
\clubpenalty=10000 

\font\afont=cmss12
\font\bfont=cmss9 
\font\cfont=cmssi9 
\font\slfont=cmssi10 
\font\hfont=cmcsc10 
\font\pfont=cmss10
\font\tfont=cmss12 at 18pt
\font\ttfont=cmss12 at 14pt
\font\abfont=cmssbx10 at 18pt
\font\abmfont=cmbxti10 at 14pt

\font\reffont=cmcsc10
\font\secfont=cmbx12
\def\subsecfont{\bf}
\def\subsubsecfont{\sl}
\font\contfont=cmssbx10 at 12pt

\font\footfont=cmr8
\def\ref#1{{\reffont#1}}

\def\sq{\quad}
\def\sqq{\qquad}
\def\cl{\centerline}

\newcount\secnum
\newcount\defnum                 
\newcount\subsecnum
\newcount\subsubsecnum
\newcount\remarknum
\newcount\lemmanum
\newcount\thmnum
\newcount\vbnum

\def\sect#1{\advance\secnum by 1\subsecnum=0 \equationnumber=0 
                     \remarknum=0\lemmanum=0\thmnum=0\defnum=0
                \vskip 12pt
                \leftline{\secfont\the\secnum.\    #1}
                \message{#1}\nobreak\noindent}
                 
	\def\subsect#1{\advance\subsecnum by 1\subsubsecnum=0
                \vskip 8pt
                \leftline{\subsecfont \the\secnum.\the\subsecnum.\ #1}
                \message{#1}\nobreak \noindent}

\def\subsubsect#1{\advance\subsubsecnum by 1
               \vskip 6pt
               \leftline{\subsubsecfont\the\secnum.\the\subsecnum.\the\subsubsecnum. \ #1}
               \message{#1}\nobreak \noindent}

\def\remark{\advance \remarknum by 1
         \vskip\baselineskip
         \noindent {\bf Remark \the\secnum.\the\remarknum.}\quad\rm
         }
\def\lemma{\advance \lemmanum by 1
         \vskip\baselineskip
         \noindent {\bf Lemma \the\secnum.\the\lemmanum.}\quad\it
         }
\def\theorem{\advance \thmnum by 1
         \vskip\baselineskip
         \noindent {\bf Theorem \the\secnum.\the\thmnum.}\quad\it
         }
\def\proof{
         \vskip\baselineskip
         \noindent {\bfProof }\quad\rm
         }
\def\vb{\advance \vbnum by 1
         \vskip\baselineskip
         \noindent {\bf Example 
\the\vbnum.}\quad
         }
\def\defi{\advance \defnum by 1
         \vskip\baselineskip
         \noindent {\bf Definition 
\the\defnum.}\quad
         }

\def\brg#1#2{{{\lower.6ex
\hbox{$\scriptstyle#1$}}\over 
{\raise.8ex
\hbox{$\scriptstyle#2$}}}}

\def\remp 
#1. #2\par{\medbreak\noindent{\bf#1.\enspace}{\rm#2}\par\medbreak} 

\def\thep
#1. #2\par{\medbreak\noindent{\bf#1.\enspace}{\sl#2}\par\medbreak} 

\def\prop
#1. #2\par{\medbreak\noindent{\bf#1.\enspace}{\rm#2}\par
 \rightline{\vrule height4pt width5.5pt depth2pt}\medbreak}
\def\proof{\bf \medbreak \noindent Proof. \rm} 
\def\eoproof{{\unskip\nobreak\hfil\penalty50
	\hskip2em\hbox{}\nobreak\hfil\vrule height4pt width5.5pt depth2pt
	\parfillskip=0pt\finalhyphendemerits=0\medbreak}}
    
\def\w#1{{\sqrt#1}\,} 
\def\kd{\partial}
\def\sq{\quad}
\def\sqs{\qquad}
\def\n{\eqno}
\def\cl{\centerline}
\def\el{\eqalign}

\def\iy{\infty}

\def\bo{{\cal O}}

\def\pd#1#2{{{\kd#1}\over{\kd#2}}}
\def\pdt#1#2{{{\kd^2#1}\over{\kd #2^2}}}

\def\br#1#2{{{#1}\over{#2}}}

\def\intp{\int_0^\iy}
\def\intr{\int_{-\iy}^\iy}

\def\RR{{{\rm I}\!{\rm R}}}
\def\NN{{{\rm I}\!{\rm N}}}
\def\RRP{{{\rm I}\!{\rm R^+}}}
\def\RRN{{{\rm I}\!{\rm R_0}}}
\def\NNP{{{\rm I}\!{\rm N^+}}}
\def\ZZ{{\hbox{Z}\!\!\hbox{Z}}}
\def\KK{{{\rm I}\!{\rm K}}}
\def\one{{{\rm 1}\hskip-0.55ex{\rm I}}}
\def\CC{\hbox{\rlap{$\,\,
  $\hbox{\vrule height6.2pt width.35pt depth-0.1pt}}$\rm C$}}
\def\QQ{\hbox{\rlap{$\,\,
  $\hbox{\vrule height6pt width.35pt depth0.1pt}}$\rm Q$}}
\def\P{\cal P}

\def\frac#1/#2{\leavevmode\kern.1em\raise.5ex\hbox{\the\scriptfont0
#1}\kern-.1em/\kern-.15em\lower.25ex\hbox{\the\scriptfont0
#2}}

%
\newcount\equationnumber	\equationnumber=0
\def\eqnum{\relax
	\global\advance\equationnumber by 1
	\equationnumberformat{\the\equationnumber}%
	}%
\def\eqname#1{\relax
	\count255=\equationnumber
	\assignnumber{EN#1}\equationnumber
	\global\equationnumber=\count255
	\global\advance\equationnumber by 1
	\ifnum\csname EN#1\endcsname=\equationnumber
	\else
		\message{The equation number for ``#1'' is incorrect!}%
	\fi
	\equationnumberformat{\csname EN#1\endcsname}%
	}%
\def\equationnumberformat#1{\eqno(\the\secnum.\equationnumbertype{#1})}%
\def\equationnumbertype#1{\number#1\relax}%
\def\referenceequation#1{\relax
	\assignnumber{EN#1}\equationnumber
	\equationnumbertype{\csname EN#1\endcsname}%
	}%
\def\forwardreferenceequation#1#2{\relax
	\global\advance\equationnumber by #2
	\assignnumber{EN#1}\equationnumber
	\global\advance\equationnumber by -1
	\global\advance\equationnumber by -#2
	\referenceequation{#1}%
	}%
%
\def\assignnumber#1#2{\relax
	\ifnum0<0\csname#1\endcsname
	\else
		\global\advance#2 by 1
		\expandafter\expandafter\expandafter
			\xdef\csname#1\endcsname{\the#2}%
	\fi
	}%

\def\fre{\forwardreferenceequation}
\def\en{\eqname}
\def\req#1{(\the\secnum.{\referenceequation{#1}})}

\def\arcsinh{{\rm arcsinh}}
\def\arccosh{{\rm arccosh}}
\def\arctanh{{\rm arctanh}}
\def\arccoth{{\rm arccoth}}

\def\erf{{\rm erf}}
\def\erfc{{\rm erfc}}

\def\phase{{\rm ph}}
\def\sign{{\rm sign}}

\def\wt{{\sqrt{2}}}

\def\phih{{\widehat \phi}}
\def\psih{{\widehat \psi}}

\def\Ai{{{\rm Ai}}}
\def\Bi{{{\rm Bi}}}

\def\phizeta{\left(\br{4\z}{1-z^2}\right)^{1/4}}

\catcode`@=11 

\font\ninerm=cmr10 at 9pt
\font\eightrm=cmr7 at 8pt
\font\sixrm=cmr7 at 6pt

\font\ninei=cmmi10 at 9pt
\font\eighti=cmmi10 at 8pt
\font\sixi=cmmi7 at 6pt
\skewchar\ninei='177 
\skewchar\eighti='177 
\skewchar\sixi='177

\font\ninesy=cmsy10 at 9pt
\font\eightsy=cmsy10 at 8pt
\font\sixsy=cmsy7 at 6pt
\skewchar\ninesy='60 
\skewchar\eightsy='60 
\skewchar\sixsy='60

\font\eightss=cmssq8

\font\eightssi=cmssqi8

\font\ninebf=cmbx10 at 9pt
\font\eightbf=cmbx10 at 8pt
\font\sixbf=cmbx7 at 6pt

\font\ninett=cmtt10 at 9pt
\font\eighttt=cmtt10 at 8pt

\hyphenchar\tentt=-1 
\hyphenchar\ninett=-1
\hyphenchar\eighttt=-1

\font\ninesl=cmsl10 at 9pt
\font\eightsl=cmsl10 at 8pt

\font\nineit=cmti10 at 9pt
\font\eightit=cmti10 at 8pt

\font\tenu=cmu10 
\newskip\ttglue

\def\ninepoint{\def\rm{\fam0\ninerm}%
  \textfont0=\ninerm \scriptfont0=\sixrm \scriptscriptfont0=\fiverm
  \textfont1=\ninei \scriptfont1=\sixi \scriptscriptfont1=\fivei
  \textfont2=\ninesy \scriptfont2=\sixsy \scriptscriptfont2=\fivesy
  \textfont3=\tenex \scriptfont3=\tenex \scriptscriptfont3=\tenex
  \def\it{\fam\itfam\nineit}%
  \textfont\itfam=\nineit
  \def\sl{\fam\slfam\ninesl}%
  \textfont\slfam=\ninesl
  \def\bf{\fam\bffam\ninebf}%
  \textfont\bffam=\ninebf \scriptfont\bffam=\sixbf
   \scriptscriptfont\bffam=\fivebf
  \def\tt{\fam\ttfam\ninett}%
  \textfont\ttfam=\ninett
  \tt \ttglue=.5em plus.25em minus.15em
  \normalbaselineskip=11pt
  \def\MF{{\manual hijk}\-{\manual lmnj}}%
  \let\sc=\sevenrm
  \let\big=\ninebig
  \setbox\strutbox=\hbox{\vrule height8pt depth3pt width\z@}%
  \normalbaselines\rm}

\def\eightpoint{\def\rm{\fam0\eightrm}%
  \textfont0=\eightrm \scriptfont0=\sixrm \scriptscriptfont0=\fiverm
  \textfont1=\eighti \scriptfont1=\sixi \scriptscriptfont1=\fivei
  \textfont2=\eightsy \scriptfont2=\sixsy \scriptscriptfont2=\fivesy
  \textfont3=\tenex \scriptfont3=\tenex \scriptscriptfont3=\tenex
  \def\it{\fam\itfam\eightit}%
  \textfont\itfam=\eightit
  \def\sl{\fam\slfam\eightsl}%
  \textfont\slfam=\eightsl
  \def\bf{\fam\bffam\eightbf}%
  \textfont\bffam=\eightbf \scriptfont\bffam=\sixbf
   \scriptscriptfont\bffam=\fivebf
  \def\tt{\fam\ttfam\eighttt}%
  \textfont\ttfam=\eighttt
  \tt \ttglue=.5em plus.25em minus.15em
\font\reffont=cmcsc10 at 9pt
  \normalbaselineskip=9pt
  \def\MF{{\manual opqr}\-{\manual stuq}}%
  \let\sc=\sixrm
  \let\big=\eightbig
  \setbox\strutbox=\hbox{\vrule height7pt depth2pt width\z@}%
  \normalbaselines\rm}

\def\a{\alpha}\def\b{\beta}\def\c{\gamma}\def\d{\delta} \def\vth{\vartheta}
\def\eps{\varepsilon}\def\f{\phi}\def\k{\kappa}\def\l{\lambda}\def\m{\mu}
\def\p{\pi}\def\r{\rho}\def\s{\sigma}\def\t{\tau}\def\th{\theta}
\def\x{\xi}\def\y{\eta}\def\z{\zeta}\def\om{\omega}   
\def\La{\Lambda}\def\oom{\Omega}\def\G{\Gamma}\def\D{\Delta}

\def\sn{\sum_{n=0}^\iy}
\def\sk{\sum_{k=0}^\iy}
\def \som#1#2{\sum_{{#1}}^{{#2}}}

\def\nfe{\bigskip}
\def\vfe{\vfil\eject}
\def\oom{\Omega}
\def\wh{\widehat}

\def\gp{$C_n^{\c}(x)$}
\def\lp{$L_n^{\a}(x)$}

\def\gpf{C_n^{\c}(x)}
\def\lpf{L_n^{\a}(x)}

\def\bri{\br1{2\pi i}}
\def\C{{\cal C}}
\def\H{{\cal H}}
\def\L{{\cal L}}

\def\index{}
\def\insil{}
\def\inrm{}

\def\ph#1{(#1)_n}

\def\be{\beta}
\def\al{\alpha}
\def\la{\lambda}

\def\aa{^{(\a)}(x)}

\def\ten{{{\,\times\,10}}}

\def\ww#1{{\eightpoint\sqrt#1}\,} 

\def\({\left(}
\def\){\right)}

\def\[{\left[}
\def\]{\right]}

\def\bxf#1{
$$
\vbox{\hrule\hbox{\vrule\kern3pt
\vbox{\kern3pt\hbox{\kern3pt$\displaystyle{#1}\kern3pt$}\kern3pt}
\vrule}\hrule}
$$
}
%
%
\def\bxfn#1#2{
$$
\el{
\vbox{\hrule\hbox{\vrule\kern3pt
\vbox{\kern3pt\hbox{\kern3pt$\displaystyle{#1}\kern3pt$}\kern3pt}
\vrule}\hrule}}
\en{#2}
$$
}

\def\F#1#2#3#4{{}_2F_1\left(\matrix{#1,#2\cr#3\cr};\,#4\right)}
\def\FFF#1#2#3#4#5#6{{}_3F_2\left(\matrix{#1,#2,#3\cr#4,#5\cr};\,#6\right)}
\def\jac#1#2#3#4{P_#1^{(#2,#3)}\left(#4\right)}
\def\lag#1#2#3{L_#1^{#2}\left(#3\right)}
\def\bxft#1{\medskip
\vbox{\hrule\hbox{\vrule\kern3pt
\vbox{\kern3pt\hbox{\kern3pt{#1}\kern3pt}\kern3pt}
\vrule}\hrule}
}

\def\llra{\Longleftrightarrow}

\bigskip 

\tfont 
\centerline         {Large Parameter Cases of the} 
\bigskip
\centerline         {Gauss Hypergeometric Function}
\bigskip\bigskip
\afont
\cl{Nico M. Temme}
\medskip 
\cfont
\cl{
CWI,  
P.O. Box 94079, 
1090 GB Amsterdam, 
The Netherlands
}
\cfont
\cl{e-mail: \tt nicot@cwi.nl}\medskip          
\centerline               {}

\bfont
\parindent=25pt 
{\pfont ABSTRACT}\par\noindent
{\narrower\noindent
We consider the asymptotic behaviour of the Gauss hypergeometric
function when several of the parameters {\it a, b, c}  are large.
We indicate which cases are of interest for orthogonal polynomials
(Jacobi, but also Krawtchouk, Meixner, etc.), which results are
already available and which cases need more attention. We also
consider a few examples of ${}_3${\it F}${}_2$ functions of unit argument,
to explain which difficulties arise in these cases, when standard
integrals or differential equations are not available.

\vskip 0.4cm \noindent 
\cfont 
2000 Mathematics Subject Classification: 
\bfont 
33C05, 33C45, 41A60, 30C15, 41A10.
\par\noindent  
\cfont 
Keywords \& Phrases: 
\bfont 
Gauss hypergeometric function,
asymptotic expansion.
\par\noindent   
{\parindent 35pt
\cfont  
\item{Note:\quad}
\bfont  Work carried out under project MAS1.2 
Analysis, Asymptotics and Computing. 
This report has been accepted for publication in the proceedings of
{\sl The Sixth International Symposium on Orthogonal Polynomials, Special
Functions and their Applications}, 18-22 June 2001, Rome, Italy.
\par
}  
}
\rm
\parindent=15pt
\sect{Introduction}%
The Gauss hypergeometric function
$$ \F abcz
=1+\br {ab}{c}z+\br{a(a+1)\,b(b+1)}{c(c+1)\,2!}z^2+\ldots
=\sn\br{\ph a\ph b}{\ph c\ n!}\,z^n,\en{i1}
$$
where
$$
\ph a
=\br{\Gamma(a+n)}{\Gamma(a)}=(-1)^n\,n!\,{-a\choose n}, \en{i2}
$$
is defined for $|z|<1$ and $c\ne 0,-1,-2,\ldots\ $.

A simplification occurs when $b=c$:
$$
\F abbz=(1-z)^{-a}.\en{i3}
$$

When $a$ or $b$ are non-positive integers the series in \req{i1}
will terminate, and $F$ reduces to a polynomial. We have 
$$\F{-n}bcz
= \som{k=0}n\br{(-n)_m\,(b)_m}{(c)_m \ m!}\,z^m
= \som{m=0}n (-1)^m{n\choose m}\br{(b)_m}{(c)_m}\,z^m. \en{i4}
$$

The value at $z=1$ is defined when $\Re(c-a-b)>0$ and is given by
$$\F abc1=\br{\G(c)\G(c-a-b)}{\G(c-a)\G(c-b)}.\en{i5}
$$

For the polynomial case we have
$$\F{-n}bc1=\br{\ph{c-b}}{\ph{c}},\sq
n=0,1,2,\ldots\ .\en{i6}
$$ 

Generalizing, let $p,q= 0,1,2,\ldots$ with $p\le q+1$. Then
$$
_pF_q\(\matrix{a_1,\ldots,a_p\cr
b_1,\ldots,b_q\cr};\;z\)=
\sn\br{\ph{a_1}\cdots\ph{a_p}}{\ph{b_1}\cdots\ph{b_q}}\,\br{z^n}{n!}.\en{i7}
$$
This series converges for all $z$ if $p<q+1$ and for $|z|<1$ if $p=q+1$.

Sometimes we know the value of a terminating  function 
at $z=1$ (Saalsch\"utz's theorem):
$${}_3F_2\(\matrix{
-n,\;a,\;b\;\cr
a,\;1+a+b-c-n\;\cr}\;;\;1\)=\br{\ph{c-a}\ph{c-b}}{\ph c\ph{c-a-b}},\en{i8}
$$
where $n=0,1,2,\ldots$\ .

The behaviour of the Gauss hypergeometric function
$F(a,b;c;z)$ for large $|z|$ follows from the transformation formula
$$\el{
\F abcz&=\br{\G(c)\G(b-a)}{\G(b)\G(c-a)}\
(-z)^{-a}\F{a}{a-c+1}{a-b+1}{\br1z} \cr
&\quad+\br{\G(c)\G(a-b)}{\G(a)\G(c-b)}\ 
(-z)^{-b}\F {b}{b-c+1}{b-a+1}{\br1{z}},\cr
}
\en{i9}$$
where $|\phase(-z)|<\pi$,
and from the expansion
$$\el{\F abcz
&=1+\br {ab}{c}z+\br{a(a+1)\,b(b+1)}{c(c+1)\,2!}z^2+\ldots\cr
&=\sn\br{\ph a\ph b}{\ph c\ n!}\,z^n,\cr
}\en{i10}
$$
with
$$
\ph a
=\br{\Gamma(a+n)}{\Gamma(a)}=(-1)^n\,n!\,{-a\choose n},\sq n=0,1,2,\ldots,
\en{i11}
$$
and where $|z|<1$. We always assume that $c\ne 0,-1,-2,\ldots\ $.

The asymptotic behaviour of $F(a,b;c;z)$ for the case that one or more of
the parameters $a$, $b$ or $c$ are large is more complicated (except when 
only $c$ is large). Several contributions in the literature are available
for the asymptotic expansions of functions of the type
$$
F(a+e_1\la,b+e_2\la;c+e_3\la;z),\quad e_j=0,\pm1,\quad \la\to\infty.\en{i12}
$$

In this paper we give an overview of these results, and we
show how this set of 26 cases can be reduced by using several types
of transformation formulas. 

In particular we indicate which cases are of interest for orthogonal 
polynomials that can be expressed in terms of the Gauss hypergeometric 
function. The Jacobi polynomial $P_n^{(\a,\b)}(x)$ 
is an important example, and we mention 
several cases of this polynomial in which $n$ is a large integer
and $\a$ and/or $\b$ are large. When $\a$ or $\b$ are
negative the zeros of $P_n^{(\a,\b)}(x)$ may be outside the 
interval  $[-1,1]$, and for several cases we give the distribution of the 
zeros in the complex plane. 

We mention a few recent papers where certain uniform 
expansions of the Gauss function are given (in terms of Bessel  functions 
and Airy functions), from which expansions the distribution of the zeros 
can be obtained. Many cases need further investigations.

Of a completely  different nature is the asymptotics of generalised 
hypergeometric functions. We consider a few cases of $_3F2$ terminating
functions of argument $1$ and $-1$, and show for these cases some ad hoc 
methods. In general, no standard methods based on integrals or 
differential equations exist for these quantities.

\sect{Asymptotics: a first example}%
Consider the asymptotics of 
$$
\F a {\b+\la}{\c+\la} z,\sq \la\to\iy. \en{s1}
$$
We use ${\b+\la}\sim{\c+\la}$ and we try (using (1.3))
$$
\F a {\b+\la}{\c+\la}z{\sim}\F a {\b+\la}{\b+\la} z=(1-z)^{-a}\ ?\en{s2}
$$
Observe first that, if $\Re (\c -a - \b)> 0$,
$$
\F a {\b+\la}{\c+\la}1=\br{\G(\c+\la)\G(\c-a-\b)}{\G(\c+\la-a)\G(\c-\b)}.\en{s3}
$$
We see that \req{s1} cannot hold for $z\sim 1$ (if $\Re (\c -a - \b)> 0$).

However, we can use
$$\F abcz =(1-z)^{-a}\F a{c-b}c{\br z{z-1}}.\en{s4}$$
This gives
$$\el{
&\F a {\b+\la}{\c+\la} z
=(1-z)^{-a}\F a{\c-\b}{\c+\la}{\br z{z-1}}\cr
&\quad=(1-z)^{-a}\Biggl[1+\br{a(\c-\b)}{\c+\la}\br z{z-1}+
\br{(a)(a+1)(\c-\b)(\c-\b+1)}{(\c+\la)(\c+\la+1)\,2!}\(\br
z{z-1}\)^2\ldots\Biggr]\cr 
&\quad=(1-z)^{-a}\[1+\br{a(\c-\b)}{\c+\la}\br z{z-1}+\bo\(\la^{-2}\)\],\cr
}
\en{s5}$$
as $\la\to\iy$ with $z$ fixed. This is the beginning of a complete asymptotic
expansion, which converges if $|z/(z-1)|<1$. It is an 
asymptotic expansion for $\la$ large, and all fixed $z, z\ne1$.

We can also use other transformation formulae:
$$\F abcz =(1-z)^{-b}\F {c-a}{b}c{\br z{z-1}}
=(1-z)^{c-a-b}\F {c-a}{c-b}c{z}.en{s6}$$
These give, with \req{s4}, the three relations
$$\el{
&\F a {\b+\la}{\c+\la} z
=(1-z)^{-a}\F a{\c-\b}{\c+\la}{\br z{z-1}}\cr
&\quad=(1-z)^{-b}\F {\c+\la-a}{\b+\la}{\c+\la}{\br z{z-1}}
=(1-z)^{\c-a-\b}\F {\c+\la-a}{\c-\b}{\c+\la}{z}\cr
}
\en{s7}$$
We see in \req{s6} that the large parameter $\la$ can be
distributed over other parameter places. In the present case only the 
second form is suitable for
giving an asymptotic expansion when using the Gauss series.
The third form gives a useless Gauss series (for large $\la$).
In the final form  the Gauss series converges at $z=1$ if $\Re(a+\b-\c)>0$,
but the series does not have an asymptotic property for large values of $\la$.

The transformation formulae in \req{s6} are an important tool for obtaining 
asymptotic expansions. We will use these and other formulae for investigating 
all cases of (1.12).
 
So far, we mention the Gauss series for obtaining an asymptotic expansion.
However,  to obtain an asymptotic expansion of \req{s1} that holds uniformly if
$z$ is close to 1, we need a different approach, as will be explained in later
sections.
\sect{Some history and recent activities}%
{\parindent=30pt
\item{[1]\sq}
\ref{Watson} (1918) studied the cases
$$
\F{\alpha+\la}{\beta+\la}{\gamma+2\la}{\br2{1-z}},\quad
\F{\alpha+\la}{\beta-\la}{c}{\br{1-z}2},\quad
\F{a}{b}{\gamma+\la}{z}.
$$
All results are of Poincar\'e type (i.e., negative powers of the large
parameter) and   hold for large domains of complex parameters and argument.
Watson began with contour integrals  and evaluated them by the method of steepest
descent.
\item{[2]\sq}
\ref{Luke} (1969, Vol. I, p. 235) summarizes these results, 
and gives many other results, also for higher $_pF_q-$functions. 
In particular Luke investigates "extended Jacobi polynomials", which are of the
type
$$
{}_{p+2}F_q\(\matrix{-n,n+\la,a_1,\ldots,a_p\cr
b_1,\ldots,b_q\cr};\;z\).
$$ 
For $p=0,q=1$ and integer $n$ this is the Jacobi polynomial.
\item{[3]\sq}
\ref{Jones} (2001) considers a uniform expansion of
$$
\F{\alpha+\la}{\beta-\la}{c}{\br12-\br12z}
$$
and gives a complete asymptotic expansion in term of $I-$Bessel functions
with error bounds. Applications are discussed for Legendre functions.
\item{[4]\sq}
\ref{Olde Daalhuis} (2001) and (2002) give new uniform expansions of 
$$
\F{a}{\beta-\la}{\gamma+\la}{-z}
$$
in terms of parabolic cylinder functions, and of 
$$
\F{\alpha+\la}{\beta+2\la}{c}{-z}
$$
in terms of Airy functions.
\item{[5]\sq}
Other recent contributions on uniform expansions for Gauss functions
follow from
\ref{Ursell} (1984) (Legendre functions),
\ref{Boyd \& Dunster} (1986) (Legendre functions), 
\ref{Temme} (1986) (see later section), 
\ref{Frenzen} (1990) (Legendre functions),
\ref{Dunster} (1991) (conical functions),
\ref{Dunster} (1999) (Jacobi and Gegenbauer polynomials),
\ref{Wong} and co-workers 
(1992) (Jacobi polynomials),
(1996) (Pollaczek),
(1997) (Jacobi function),
(1998) and (1999) (Meixner), 
(2000) (Kraw\-tchouk),
(2001) (Meixner-Pollaczek).

\item{[6]\sq}
On Bessel polynomials, which are of $_2F_0-$type, see 
\ref{Wong \& Zhang} (1997) and \ref{Dunster} (2001).
On Charlier polynomials ($_2F_0-$type) see 
\ref{Bo Rui \& Wong} (1994). On Laguerre polynomials
($_1F_1-$type) see \ref{Frenzen \& Wong} (1988) 
and \ref{Temme} (1990).
\par
}

%
%
$$\vbox{\hsize=13.4truecm
\noindent
$$\vbox{\offinterlineskip
\halign{\strut        \vrule
\hfil $\ #$\hfil &    \vrule 
\hfil $\ #$ \hfil &   \vrule
\hfil $\ #$ \hfil &    \vrule 
${#\ }$          &   \vrule
\hfil $\ #$\hfil &    \vrule 
\hfil $\ #$ \hfil &   \vrule
\hfil $\ #$ \hfil  &   \vrule 
${#\ }$         &    \vrule
\hfil $\ #$\hfil &    \vrule 
\hfil $\ #$ \hfil &   \vrule
\hfil $\ #$ \hfil     \vrule 
\cr\noalign{\hrule}
\hfil \quad e_1 \quad\hfil&\hfil\quad e_2\quad\hfil&\hfil \quad
e_3\quad\hfil&&
\hfil \quad e_1 \quad\hfil&\hfil\quad e_2\quad\hfil&\hfil \quad
e_3\quad\hfil&&
\hfil \quad e_1 \quad\hfil&\hfil\quad e_2\quad\hfil&\hfil \quad
e_3\quad\hfil
\cr
\noalign{\hrule}
\    &\    &\  &\    &\    &\ & \    &\    &\ & \    &\     \cr
0  &  0  &  0   & & +  &  0  &  0 & & -  &  0  &  0 \cr 
0  &  0  &  +   & & +  &  0  &  + & & -  &  0  &  + \cr 
0  &  0  &  -   & & +  &  0  &  - & & -  &  0  &  - \cr 
0  &  +  &  0   & & +  &  +  &  0 & & -  &  +  &  0 \cr 
0  &  +  &  +   & & +  &  +  &  + & & -  &  +  &  + \cr 
0  &  +  &  -   & & +  &  +  &  - & & -  &  +  &  - \cr 
0  &  -  &  0   & & +  &  -  &  0 & & -  &  -  &  0 \cr 
0  &  -  &  +   & & +  &  -  &  + & & -  &  -  &  + \cr 
0  &  -  &  -   & & +  &  -  &  - & & -  &  -  &  - \cr 
\noalign{\hrule}
}} $$ 
{{\bf Table 1.\quad}\ All 27 cases of 
$_2F_1(a+e_1\la,b+e_2\la;c+e_3\la;z),$ \ $e_j=0,\pm1 \ .$
In the table we write $e_j=0,\pm\ .$}
}$$

%
\sect{From 27 cases to only a few}
Skipping the dummy $(\ 0 \quad 0 \sq 0\ )$, using the symmetry between the
$a$ and
$b$ parameters:
$$
\F abcz =\F bacz, \en{a1}
$$
and by using the transformation formulae of (2.6)
we reduce the 27 cases to the 8 cases:
$$\vbox{\hsize=9.0truecm
\noindent
$$\vbox{\offinterlineskip
\halign{\strut\vrule
\hfil $\ #\quad$\hfil   &      \vrule 
\hfil $\ #$\hfil   &      \vrule 
\hfil $\ #$ \hfil   &     \vrule
\hfil $\ #$ \hfil     \vrule 
    \cr \noalign{\hrule}
&\hfil \quad e_1 \quad\hfil  &  \hfil\quad e_2\quad\hfil  &  \hfil \quad
e_3\quad\hfil    \cr 
\noalign{\hrule}
\   &   &  \      &  \     \cr 
1 & 0  &  0  &  +    \cr 
2 & 0  &  0  &  -    \cr 
3 & 0  &  +  &  0    \cr 
4 & 0  &  +  &  -    \cr 
5 & 0  &  -  &  +    \cr 
6 & +  &  +  &  -    \cr 
7 & +  &  -  &  0    \cr 
8 & -  &  -  &  +    \cr 
\noalign{\hrule}
}} $$ 
}$$

This set is obtained by using \req{a1} and (2.6), that is, only 4 of
Kummer's 24 solutions of the hypergeometric differential equation.

A further reduction of the set of 8 cases can be obtained by other solutions
of the differential equation.

Using the connection formula
$$
\el{\F abcz&=
-\br{\G(c-1)\G(a-c+1)\G(b-c+1)}{\G(a)\G(b)\G(1-c)}
  z^{1-c}(1-z)^{c-a-b}\F{1-a}{1-b}{2-c}z \cr
\quad&+\br{\G(b-c+1)\G(a-c+1)}{\G(a+b-c+1)\G(1-c)}
\F ab{a+b-c+1}{1-z},\cr
}\en{a3}
$$
we see that the second case $(\ 0 \quad 0 \sq - \ )$
in fact reduces to $(\ 0 \quad 0 \sq +\ )$.

Similarly, we have
$$
\el{
\F abcz&=
e^{a\pi i}\br{\G(c)\G(b-c+1)}{\G(a+b-c+1)\G(c-a)}
   z^{-a}\F{a}{a-c+1}{a+b-c+1}{1-\br1z} \cr
\quad&+\br{\G(c)\G(b-c+1)}{\G(a)\G(b-a+1)}
z^{a-c}(1-z)^{c-a-b}\F {1-a}{c-a}{b-a+1}{\br1{z}},\cr
}\en{a4}
$$
which reduces the third case  $(\ 0 \quad + \sq 0 \ )$
again to $(\ 0 \quad 0 \sq +\ )$.

Finally, we use
$$
\el{
\F abcz&=
e^{c\pi i}\br{\G(c-1)\G(b-c+1)\G(1-a)}
{\G(b)\G(c-a)\G(1-c)}\
z^{1-c}(1-z)^{c-a-b}\F{1-a}{1-b}{2-c}z \cr
&+e^{(c-a)\pi i}\br{\G(1-a)\G(b-c+1)}{\G(1-c)\G(b-a+1)}
 z^{a-c}(1-z)^{c-a-b}\F {1-a}{c-a}{b-a+1}{\br1{z}},\cr
}\en{a5}
$$
which reduces the forth case  $(\ 0 \quad + \sq - \ )$
to the fifth case $(\ 0 \quad - \sq +\ )$.

Consequently we have the five remaining cases:
$$\vbox{\hsize=9.0truecm
\noindent
$$\vbox{\offinterlineskip
\halign{\strut\vrule
\hfil $\ #\quad$\hfil   &      \vrule 
\hfil $\ #$\hfil   &      \vrule 
\hfil $\ #$ \hfil   &     \vrule
\hfil $\ #$ \hfil     \vrule 
    \cr \noalign{\hrule}
&\hfil \quad e_1 \quad\hfil  &  \hfil\quad e_2\quad\hfil  &  \hfil \quad
e_3\quad\hfil    \cr 
\noalign{\hrule}
\   &   &  \      &  \     \cr 
1 & 0  &  0  &  +    \cr 
5 & 0  &  -  &  +    \cr 
6 & +  &  +  &  -    \cr 
7 & +  &  -  &  0    \cr 
8 & -  &  -  &  +    \cr 
\noalign{\hrule}
}} $$ 
}$$

These five cannot be reduced further by using relations between the functions
$$
\F{a+e_1\la}{b+e_2\la}{c+e_3\la}z,\sq
e_j=0,\pm1.$$

However (\ref{Olde Daalhuis} (2001)), using 
$$
\el{
\F abcz&=
e^{a\pi i}\br{\G(c)\G(b-a)}{\G(b)\G(c-a)}
z^{-a}\F{a}{a-c+1}{a-b+1}{\br1z} \cr
&+e^{b\pi i}\br{\G(c)\G(a-b)}{\G(a)\G(c-b)}
 z^{-b}\F {b}{b-c+1}{b-a+1}{\br1{z}},\cr
}\en{a6}
$$
and 
$$
\el{
\F abcz&=
e^{(c-b)\pi i}\br{\G(c)\G(b-a)}{\G(b)\G(c-a)}
z^{b-c}(1-z)^{c-a-b}\F{1-b}{c-b}{a-b+1}{\br1z} \cr
&+e^{(c-a)\pi i}\br{\G(c)\G(a-b)}{\G(a)\G(c-b)}
 z^{a-c}(1-z)^{c-a-b}\F {1-a}{c-a}{b-a+1}{\br1{z}},\cr
}\en{a7}
$$
we see that the cases $(\ + \quad + \sq - \ )$ and 
$(\ - \quad - \sq + \ )$ both can be handled by the case
$(\ + \quad 2+ \sq 0 \ )$.

Hence, the smallest set is the following set of four:
$$\vbox{\hsize=9.0truecm
\noindent
$$\vbox{\offinterlineskip
\halign{\strut\vrule
\hfil $\ #\quad$\hfil   &      \vrule 
\hfil $\ #$\hfil   &      \vrule 
\hfil $\ #$ \hfil   &     \vrule
\hfil $\ #$ \hfil     \vrule 
    \cr \noalign{\hrule}
&\hfil \quad e_1 \quad\hfil  &  \hfil\quad e_2\quad\hfil  &  \hfil \quad
e_3\quad\hfil    \cr 
\noalign{\hrule}
\   &   &  \      &  \     \cr 
A & 0  &  0  &  +    \cr 
B & 0  &  -  &  +    \cr 
C & +  &  -  &  0    \cr 
D & +  &  2+  &  0    \cr 
\noalign{\hrule}
}} $$ 
}$$
if we accept the special case $D$:
$$\F{a+\la}{b+2\la}{c}z.$$

\subsect{Quadratic transformations}
For special combinations of the parameters $a, b, c$ other types of relations 
between $_2F_1-$ functions exist. For example, we have the quadratic
transformation
$$
\F{a}{b}{2b}{z}=\(1-z\)^{-a/2}\F{\brg12a}{b-\brg12a}{b+\brg12}
{\br{z^2}{4z-4}}.
$$
For large $a$ the left-hand side is of type 
$$
(\ + \ 0 \ 0 \ ) \equiv (\ 0 \ + \ 0 \ ) \equiv (\ 0 \ 0 \ + \ ) = A,
$$
whereas the right-hand side is of type $C = (\ + \ - \ 0 \ )$.
\sect{Which special functions are involved?}
The next step is to indicate which special functions are associated 
with the above four cases. 
\subsect{Legendre functions}%
We have the Legendre function
$$
P_\nu^{\mu}(z)=\br{(z+1)^{\br12\mu}(z-1)^{-\br12\mu}}{\G(1-\mu)}
\F{-\nu}{\nu+1}{1-\mu}{\br12(1-z)}. \en{s1}
$$

Case $A:\quad (\ 0 \ 0 \ + \ )$

We see that case $A$ occurs 
when $\mu\to-\iy$. Several relations between 
the Legendre function and the Gauss function give also the case 
$\mu\to+\iy$.

Case $B:\quad (\ 0 \ - \ + \ )$

Because 
$$(\ - \ + \ -\ )  \equiv  ( \ 0 \ + \ -\ )  
\equiv  ( \ 0 \ - \ +\ )  =  B $$
(see (3.5)), we see that case $B$ occurs when 
in $P_\nu^{\mu}(z)$ both parameters 
$\nu$ and $\mu$ tend to $+\iy$.

Case $C:\quad (\ + \ - \ 0\ ) \equiv (\ - \ + \ 0\ ) $

This also occurs for 
$P_\nu^{\mu}(z)$, with $\nu$ large  and $\mu$  fixed. 

In the results for Legendre functions in the literature 
rather flexible conditions on the parameters 
$\mu$ and $\nu$ are allowed. For example, in \ref{Dunster} (1991)
it is assumed that the ratio $\nu/\mu$ is bounded.

Any hypergeometric function, for which a quadratic transformation exists,
can be expressed in terms of Legendre functions. So, many special cases, and
mixed versions of the cases $A, B, C, D$ are possible for Legendre functions.  
\subsect{Hypergeometric polynomials}%
Several  orthogonal polynomials of hypergeometric type
(see the Askey scheme in \ref{Koekoek \& Swartouw} (1998))
have representations 
with Gauss hypergeometric functions. Next to the special cases of 
the Jacobi polynomials (Gegenbauer, Legendre) we have 
the following cases.

\bxft{{\bf Pollaczek}}
$$P_n(\cos\theta;a,b)=e^{n\theta i}
\F{-n}{\br12+i\phi(\theta)}{1}{1-e^{-2i\theta}},$$
$$ \phi(\theta)=\br{a\cos\theta + b}{\sin\theta}.$$
For uniform expansions in terms of Airy functions see \ref{Bo Rui \& Wong} (1996).
The asymptotics is for $n\to\iy, \theta=t/\sqrt{n}$ with $t$ bounded, and bounded away
from zero. 

\bxft{{\bf Meixner-Pollaczek}}
$$P_n^{(\la)}(x;\phi)=\br{(2\la)_n}{n!} e^{n\phi i}\F{-n}{\la+ix}{2\la}{1-e^{2i\phi}}.$$
For uniform expansions in terms of parabolic cylinder functions
see \ref{Li \& Wong} (2001). The asymptotics is for $n\to\iy, x=\al n$ 
with $\al$ in a compact interval containing the point $\al=0$. 

\bxft{{\bf Meixner}}
$$M_n(x;\b,c)=\F{-n}{-x}{\b}{1-\br1c}.$$
In \ref{Jin \& Wong} (1998) two uniform expansions in terms of parabolic cylinder
functions are given. The asymptotics is for $n\to\iy, x=\al n$ with $\al>0$
in compact intervals containing the point $\al=1$. 

\bxft{{\bf Krawtchouk}}
$$K_n(x;p,N)=\F{-n}{-x}{-N}{\br1p},\sq n=0,1,2,\ldots,N.$$
A uniform expansion in terms of parabolic cylinder functions is given in 
\ref{Li \& Wong} (2000). The asymptotics is for $n\to\iy, x=\la N$ with 
$\la$ and $\nu=n/N$ in compact intervals of $(0,1)$.

\subsect{Jacobi polynomials}%
The Jacobi polynomial has the representation:
$$
P_n^{(\a,\b)}(x)={n+\a\choose n}\F{-n}{\a+\b+n+1}{\a+1}{\br12(1-x)}, \en{s2}
$$
and we see that case 
$C:\quad (\ + \ - \ 0\ ) \equiv (\ - \ + \ 0\ ) $ applies
if $n\to\iy$. Uniform expansions use Bessel functions for describing the asymptotic
behaviour at the points $x=\pm 1$ (Hilb-type formula, see \ref{Szeg\"o} (1975) or
for complete asymptotic expansions \ref {Wong \& Zhang} (1996)).

A well-known limit is
$$
\lim_{\al\to\iy}\al^{-n/2} P_n^{(a+\al,b+\al)}\(\br{x}{\sqrt{\al}}\)=\br{H_n(x)}{2^n\,n!},
$$
where $H_n(x)$ is the Hermite polynomial, a special case of the parabolic
cylinder function. Approximations of $H_n(x)$  are available in terms
of Airy functions. 

Several other limits are known (see \ref{Koekoek \& Swartouw} (1998)).
In \ref{Chen \& Ismail} (1991) and \ref{Dunster} (1999) asymptotic expansions 
are given for large positive parameters  $\al$ and/or $\beta$.

Of particular interest are asymptotic expansions for large $n$ and large
negative values of the parameters $\al, \be$. See 
\ref{Bosbach \& Gawronski} (1999) and 
\ref{Mart\'\i nez-Finkelshtein et al.} (2000).

With $\a=\b=-n$ we have in \req{s2} the case 
$$ (\ - \ - \ -\ ) \equiv (\ 0 \ 0 \ -\ ) \equiv (\ 0 \ 0 \ +\ ) = A.$$

These are  non-classical values for the parameters. We have, taking 
$\a=\b=-\brg12-n$, 
$$
P_n^{(\a,\b)}(z)=\br{n!}{2^{2n}}C_n^{-n}(x)=(-1)^n\sum_{m=0}^{\lfloor n/2\rfloor}
\br{n!n!2^{-n-2m}x^{n-2m}}{m!m!(n-2m)!}.
$$
This representation shows that all zeros are on the imaginary axis.

Other
representations can be used:
$$
\el{
&P_n^{(\a,\b)}(x)=
{n+\a\choose n}\(\br{1+x}{2}\)^{-\b}
\F{\a+n+1}{-\b-n}{\a+1}{\br12(1-x)}\cr
&\qquad\quad\quad={n+\b\choose n}\(\br{1-x}{2}\)^{-\a}(-1)^n
\F{\b+n+1}{-\a-n}{\b+1}{\br12(1+x)},\cr
}\en{s3}
$$
and we see that case 
$D:\ (\ + \ 2+\ \ 0\ )$ applies if $n\to\iy$ and, in the first case, $\a$ is constant
and  $\b$ equals $b-3n$, and similarly in the second case.

See Figure 1 for
the zero distribution of the Jacobi polynomial
$$
P_n^{(\a,\b)}(z), \sq n=30, \sq \a=\brg12,\sq \b=-3n+1.
$$
\centinsert{\pscaption{\psannotate{\psboxto(0.9\hsize;0cm){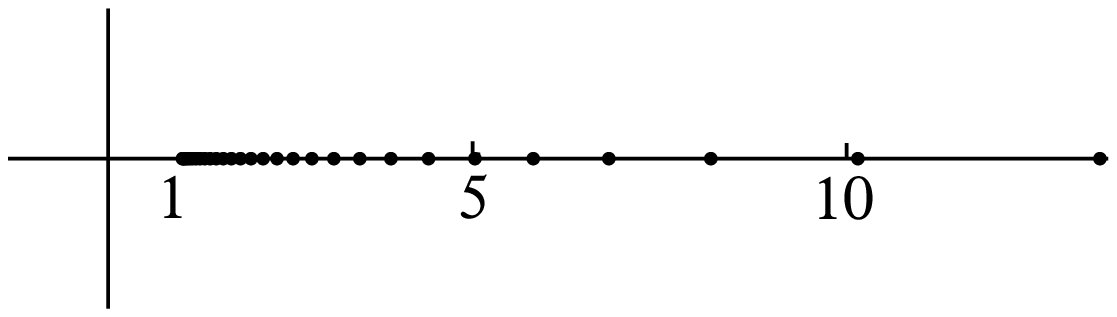
}}
{
}}{%
{\bf Figure 1.}\sq  
Zeros of 
$P_n^{(\a,\b)}(z), \ n=30, \a=\brg12,\ \b=-3n+1.$ 
}} 


The case $B:\ (\ 0 \ -\ \ +\ ) \equiv (\ - \ 0\ \ +\ ) $ applies in the representations
$$
\el{
&P_n^{(\a,\b)}(x)=\
{\a+\b+2n\choose n}\(\br{x-1}{2}\)^{n}
\F{-n}{-\a-n}{-2n-\a-\b}{\br2{1-x}}\cr
&\qquad\quad\quad={\a+\b+2n\choose n}\(\br{x+1}{2}\)^{n}
\quad\F{-n}{-\b-n}{-2n-\a-\b}{\br2{1+x}},\cr
}\en{s4}
$$
when $\a=-a-n, \b=-b-2n$ (first case) or $\b=-b-n, \a=-a-2n$ (second case), with, again,
non-classical parameter values. \ref{Olde Daalhuis} (2002) gives an expansion of the
Gauss function for this case in terms of parabolic cylinder functions. 

See Figure 2 for
the zero distribution of the Jacobi polynomial
$$
P_n^{(\a,\b)}(z), \sq n=25, \sq \a=-n+\brg12,\sq \b=-2n+1.
$$

As $n\to \iy$, the zeros approach the curve defined by
$\left|1-\(\br{3-z}{1+z}\)^2\right| = 1$. 

\centinsert{\pscaption{\psannotate{\psboxto(0.5\hsize;0cm){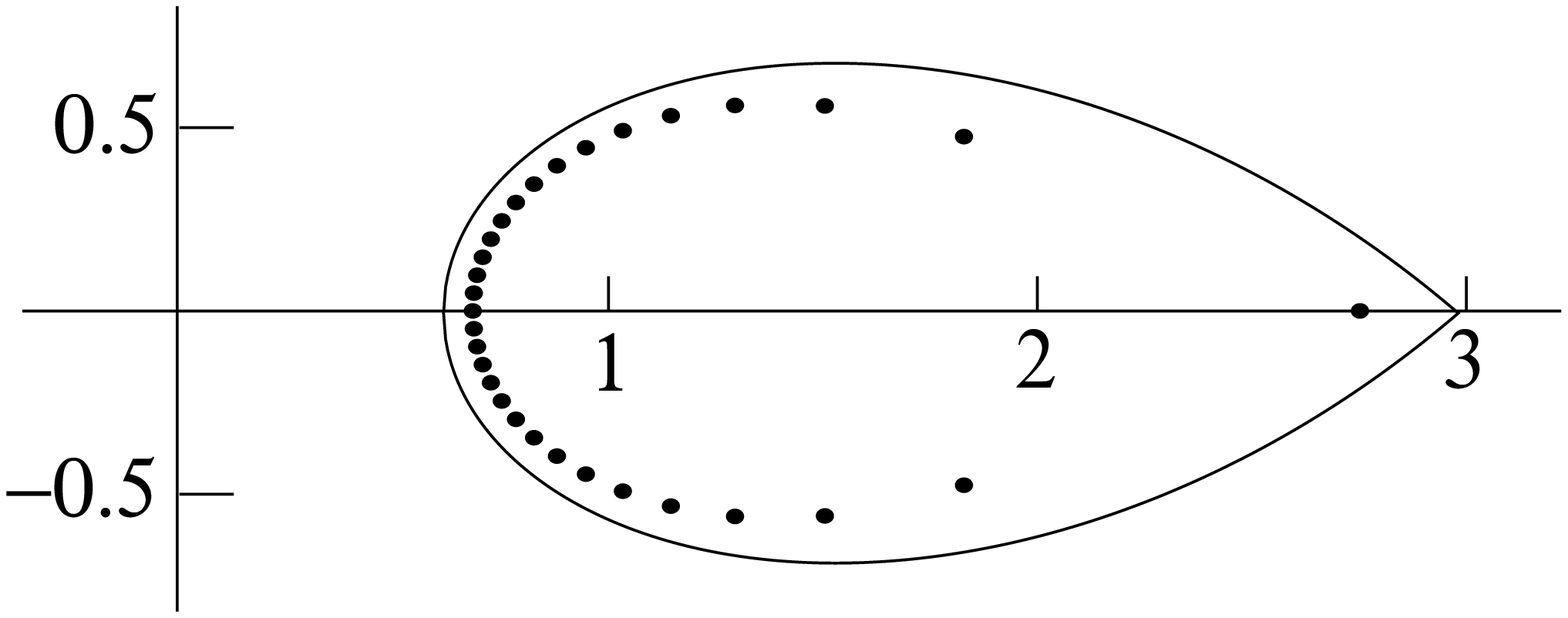
}}
{
}}{%
{\bf Figure 2.}\ 
Zeros of 
$P_n^{(\a,\b)}(z), \ n=25, \a=-n+\brg12,\ \b=-2n+1.$ 
}}

Uniform expansions for Jacobi polynomials are available in terms
of $J-$Bessel functions. Details on uniform expansions of  Jacobi polynomials in
terms of Bessel functions can be found in \ref{Wong \& Zhang} (1996),
\ref{Dunster} (1999), and in cited references in these papers.
\subsect{Gegenbauer polynomials}%
For Gegenbauer (ultraspherical) polynomials we have
several representations. For example,
$$
C_n^{\la}(x)=\br{(2\la)_n}{n!}\F{-n}{n+2\la}{\la+\brg12}{\br{1-x}{2}},
$$
and for $\la=-n$, $n\to\iy$, we have
$$
(\ - \ - \ - \ ) \equiv (\ 0 \ 0 \ - \ ) \equiv (\ 0 \ 0 \ + \ ) = A.
$$
We also have
$$
C_{2n}^{\la}(x)=(-1)^n\br{(\la)_n}{n!}\F{-n}{n+\la}{\brg12}{x^2},
$$
and for $\la=-2n$, $n\to\iy$, we have
$$
(\ - \ - \ 0\ ) \equiv (\ + \ - \ 0 \ ) = C.
$$
Also in this case a "sort of$\,$" quadratic transformation is used.

\subsect{Asymptotics of the $A-$case}%
We give a few steps in deriving the expansion for the case 
$A = (\ 0 \ 0 \ + \ )$. As always, we can use the differential equation or 
an integral representation. We take
$$
\F ab{c}z=\br{\G(c)}{\G(b)\G(c-b)}\int_0^1
\br{t^{b-1}(1-t)^{c-b-1}}{(1-zt)^{a}}\,dt,\en{c1}
$$
where $\Re c>\Re b>0, |\phase (1-z)| <\pi$. A simple transformation
gives the Laplace transform representation
$$
\F ab{c+\la}z=\br{\G(c+\la)}{\G(b)\G(c+\la-b)}\intp
t^{b-1}e^{-\la t} f(t)\,dt,
$$
with $\la\to\iy$ and where
$$
f(t)=\(\br{e^t-1}{t}\)^{b-1} e^{(1-c)t}\(1-z+ze^{-t}\)^{-a}.
$$

This standard form in asymptotics can be expanded by using Watson's lemma
(see \ref{Olver} (1974) or \ref{Wong} (1989)). For fixed $a, b, c$ and
$z$:
$$
\F ab{c+\la}z\sim\br{\G(c+\la)}{\G(c+\la-b)}
\sum_{s=0}^\iy f_s(z) \br{(b)_s}{\la^{b+s}},
$$
where $f_0(z)=1$ and other coefficients follow from 
$$f(t)=\sum_{s=0}^{\iy}
f_s(z)t^s,\sq |t|<\min(2\pi, |t_0|), \sq t_0=\ln\br{z}{z-1}.
$$ 

When $z$ is large, the singularity $t_0$ is close to the origin, and the
expansion becomes useless.

We can write
$$
\F ab{c+\la}z=\br{\G(c+\la)}{\G(b)\G(c+\la-b)}\intp
\br{t^{b-1}e^{-\la t} g(t)}{(t+\zeta)^{a}}\,dt,
$$
where 
$$\z=-t_0=\ln\br{z-1}{z}, \sq g(t)=(t+\z)^a f(t),$$
with $g$ regular at $t=t_0$. 
Expanding $g(t)=\sum g_s(z) t^s$ gives
$$
\F ab{c+\la}z \sim\br{\G(c+\la)}{\G(c+\la-b)}\ 
\sum_{s=0}^\iy g_s(z) (b)_s\z^{b-a+s}U(b+s, b-a+1+s,\z\la),
$$
in which $U$ is the confluent hypergeometric function. This expansion holds if
$\la\to\iy$, uniformly with respect to small values of $\z$, that is, large
values of $z$.

This method can be used to obtain an expansion of the Legendre function
$P_\nu^{\pm m}(z)$ for large values of $m$; see \ref{Gil, Segura \& Temme} (2000),
where we used this type of expansion for computing this Legendre function for
large $m$ and $z$. 

Other sources with related expansions for this type of Legendre functions
are \ref{Boyd \& Dunster} (1986), \ref{Dunstger} (1991), \ref{Frenzen} (1990),
\ref{Olver} (1974) and \ref{Ursell} (1984).

\subsect{Another $A-$case}%
The case $(\ - \ 0 \ 0 \ ) $ is very important for all kinds of orthogonal
polynomials and can be reduced to the previous case $A = (\ 0 \ 0 \ + \ )$.
It is of interest to give a direct approach of
$$
\F{-n}bc{z},\sq n\to\iy,\en{c2}$$
whether or not $n$ is an integer. 
First we recall the limit
$$
 \lim_{b\to\iy} \F abc{z/b}={}_1F_1\(\matrix{a\cr c\cr}\;;\;z\)
$$
which may be used as a definition for the ${}_1F_1-$function 
(the Kummer or confluent hypergeometric function). 
We are interested in the asymptotics behind this limit, 
and we expect a role of the Kummer function when in \req{c2} $n$ becomes large, 
and $z$ is small.
We again can use \req{c1}, and consider
$$
I_n=\int_0^1
{t^{b-1}(1-t)^{c-b-1}}{(1-tz)^{n}}\,dt,
$$
where $\Re c>\Re b>0, |\phase(1-z)| <\pi$. We transform 
$$1-t z=(1-z)^u=e^{u \ln(1-z)}.$$
Then we have
$$
\el{
I_n&= (1-z)^{c-b-1}\[\br{\ln(1-z)}{-z}\]^{c-1}J_n,\cr
J_n&= \int_0^1 f(u) u^{b-1} (1-u)^{c-b-1} e^{\omega u} \,du,\cr
f(u)&=\[\br{1-(1-z)^u}{-u\ln(1-z)}\]^{b-1}\ \[\br{1-(1-z)^{u-1}}{(1-u)\ln(1-z)}\]^{c-b-1},\cr
\omega&= (n+1)\ln(1-z).\cr}
$$
The function $f(u)$ is holomorphic in a neighborhood of $[0,1]$. Singularities
occur at
$$ u_k=\br{2k\pi i}{\ln(1-z)},\sq v_m=1+\br{2m\pi i}{\ln(1-z)}, \sq
k,m\in\ZZ\setminus\{1\}.
$$
So, when $z$ ranges through compacta of $\CC\setminus\{1\}$, the singularities
of $f$ are bounded away from $[0,1]$.
If $\Re\omega <0$ the dominant point in the integral is $u=0$; if $\Re\omega>0$,
then $u=1$ is the dominant point.
To obtain an asymptotic expansion for large $n$ that combines both cases,
and which will give a uniform expansion for all $z$, $|z-1|>\delta>0$, 
$f$ should be expanded at both end-points $0$ and $1$.

More details can be found in \ref{Temme} (1986), where the results 
have been applied to
a class of polynomials biorthogonal on the unit circle.

\subsect{Asymptotics of the $B-$case}%
We consider a simple case:
$$
\F {-n}{1}{n+2}{-z}=(n+1) \int_0^1 \[(1-t)(1+zt)\]^n\,dt,\en{c4}
$$
with $z$ near the point $1$.

We write 
$$\[(1-t)(1+zt)\]^n=e^{-n \phi(t)},\sq \phi(t)=-\ln\[(1-t)(1+zt)\],
$$
and we have
$$
\phi'(t)=\br{-2tz-1+z}{(1-t)(1+zt)}.
$$
We see that the integrand has a peak value at $t_0=\br{z-1}{2z}$.

So, if $z=1$ the peak is at $t=0$, if $z>1$ then $t_0\in(0,1)$,
and if $z<1$ then $t_0<0$. See Figure 3.

\centinsert{\pscaption{\psannotate{\psboxto(0.8\hsize;0cm){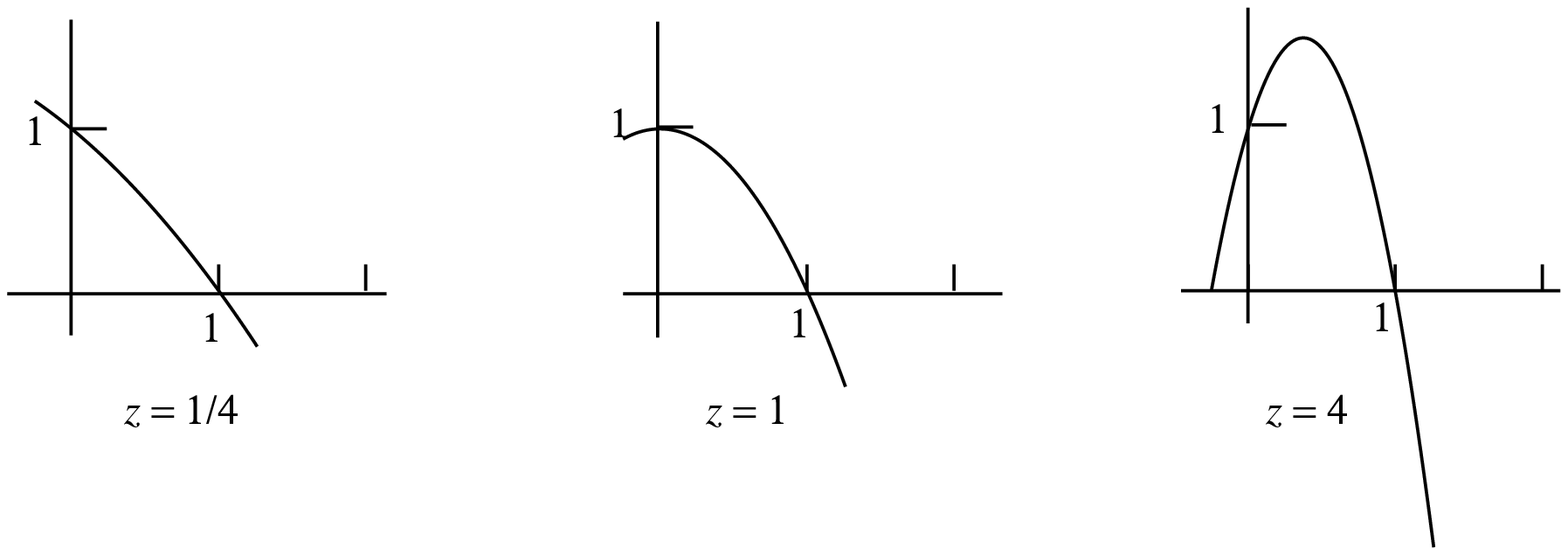
}}
{
}}{%
{\bf Figure 3.}\sq  
 $(1-t)(1+zt)$ for a few values of $z$. 
}} 

The same situation occurs for the integral
$$
\intp e^{-n(\br12w^2-\al w)}\,dw \en{c5}
$$
which has  a peak value at $w_0=\al$. This integral is an error function 
(a special case of the parabolic cylinder functions). 

In a uniform expansion of \req{c4} that holds for large $n$ and $z$
in a neighborhood of $1$ we need an error function. This explains 
why in Case $B$ (in a more general form than \req{c4}) parabolic cylinder
functions may occur for special values of $z$.

We transform the integral in \req{c4} by writing
$$
(1-t)(1+zt)=e^{-\br12w^2+\al w}
$$
with the conditions
$$
t=0 \llra w=0, \sq t=1 \llra w = \iy, \sq t=t_0 \llra w=\al.
$$
The quantity $\al$ follows from satisfying the matching of $t_0$ with $\al$:
$$
(1-t_0)(1+zt_0)=e^{-\br12w_0^2+\al w_0},
$$
giving
$$
\brg12\al^2=-\ln\[1-\(\br{z-1}{z+1}\)^2\],\sq \sign(\al)=\sign(z-1).
$$

We obtain
$$
\F {-n}{1}{n+2}{-z}=(n+1) \intp  e^{-n(\br12w^2-\al w)} f(w)\,dw,
\en{c6}
$$
where 
$$f(w)=\br{dt}{dw}=\br{w-\al}{t-t_0}\,\br{(1-t)(1+zt)}{2z}.$$

A first approximation follows by replacing $f(w)$ by 
$$f(\al)=\br{1+z}{2\sqrt{2}\,z},$$ 
giving
$$
\F {-n}{1}{n+2}{-z}\sim n f(\al) \intp  e^{-n(\br12w^2-\al w)} \,dw
=\sqrt{\pi n} \br{1+z}{4z} e^{-\br12n\al^2}
\erfc\(-\al\sqrt{n/2}\).
$$
as $n\to\iy$, uniformly with respect to $z$ in a neighborhood of $z=1$.

This expansion is in agreement with the general case 
$
_2F_1({a}, {b-\la};{c+\la};{-z})
$
considered  in \ref{Olde Daalhuis} (2002).

The relation with the Jacobi polynomials is:

$$
\F {-n}{1}{n+2}{-z}=\br{(n+1)!}{2^{2n}\,\ph{\br32}} (1+z)^n
P_n^{(n+1,-n-1)}\(\br{1-z}{1+z}\).$$
See Figure 4 for the distribution of the zeros of
$$
P_n^{(\al,\be)}(z),\sq n=30,\sq \al =n+1,\sq \be = -n-1.
$$
As $n\to \iy$, the zeros approach the curve 
$\left|1-z^2\right| = 1$.

\centinsert{\pscaption{\psannotate{\psboxto(0.5\hsize;0cm){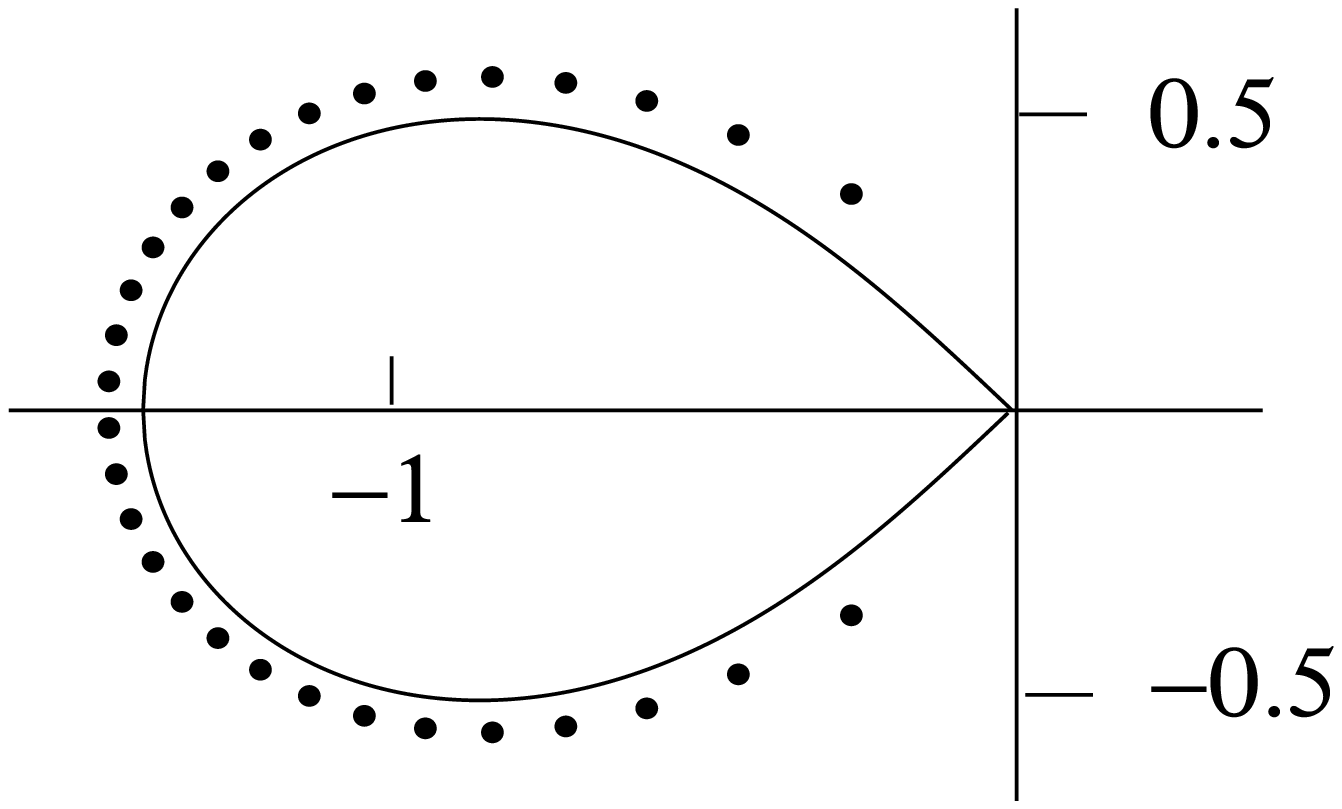
}}
{
}}{%
{\bf Figure 4.}\sq  
Zeros of 
$P_n^{(\a,\b)}(z), \ n=30, \a=-\b=n+1.$ 
}} 
\sect{Asymptotics of some ${}_3F_2$ polynomials}
In the previous examples we used integral representations of the Gauss
function. The differential equation can also be used for obtaining asymptotic
expansions. As a rule, methods based on differential equations provide more
information on the coefficients and remainders in the expansions than methods
based on integral representations. 

For the generalized hypergeometric functions asymptotic expansion are known when
the argument $z$ is large, and the other parameters are fixed; see \ref{Braaksma} (1963)
and \ref{Luke} (1969).

For the ${}_pF_q-$functions a differential equation is available of order
$\max(p,q+1)$. However, for higher order equations no methods are available for
deriving (uniform) expansions for large values of parameters.

The ${}_pF_q-$functions can be written as a Mellin-Barnes contour integral, but,
again, no methods are available for
deriving (uniform) expansions for large values of parameters from these integrals.

The terminating
${}_pF_q-$functions (one of the parameters
$a_k$ is equal to a non-positive integer) and of unit argument $(z=1)$ are of
great interest in special function theory and in applications.

In many cases recursion relations are available. For the examples to be
considered here we use ad hoc methods. It is of interest to investigate 
recursion methods for the same and other problems.

\subsect{A first example}%
As a first example (\ref{Larcombe \& Koornwinder} (2000)), consider the
asymptotics of 
$$f(n)=
\FFF{-n}{\brg12}{\brg12}{\brg12-n}{\brg12-n}{-1}
$$
It was conjectured (by Larcombe) that
$\lim_{n\to\iy} f(n)=2$
and Koornwinder gave a proof.
By using
$$
\br{\(\brg12\)_k}{\(\brg12-n\)_k}=
\br{(-1)^kn!}{\w{{\pi}}\G(n+\brg12)}\int_0^1t^{k-\br12}(1-t)^{n-k-\br12}\,dt,
$$
$k=0,1,\ldots, n$,
we can write $f(n)$ as an integral
$$f(n)=\br{n!}{\pi\,\ph{\brg12}}\int_0^1t^{-\br12}(1-t)^{n-\br12}
\F{-n}{\brg12}{\brg12-n}{\br{t}{1-t}}\,dt,
$$
and the Gauss function is a Legendre polynomial. We obtain
$$f(n)=\br{2^{-n}n!n!}{\pi\,\ph{\brg12}\,\ph{\brg12}}\int_0^\pi
\sin^n\theta\  P_n\(\br1{\sin\theta}\)\,d\theta.\en{lf1}$$
There is no straightforward way to obtain asymptotics out of this.  
A few manipulations with \req{lf1} give the result
$$
f(n)=\br{n!\,n!\,n!}{2^n\,\ph{\brg12}\,\ph{\brg12}} \ c_n$$
where $c_n$ is the coefficient in the Bessel function expansion
$$
\[e^{w/2}\,I_0(w/2)\]^2 = \sn c_n w^n,$$
and now the asymptotics easily follows. We have
$$c_n=\br1{2\pi i}\int_{\C}\br{\[e^{w/2}\,I_0(w/2)\]^2}{w^{n+1}}\, dw
=\br1{2\pi i}\int_{\C}\br{e^{2w}}{w^{n+1}}\,f(w)\, dw,\sq 
f(w)=\[e^{-w/2}\,I_0(w/2)\]^2$$
and the contour $\C$ is a circle around the
origin. The main contribution comes from the saddle point of 
$\br{e^{2w}}{w^{n+1}}$, that is from $w=w_0=n/2$.

Because $f(w_0)\sim \br2{\pi n}$, we have 
$$c_n\sim\br2{\pi n}\br{2^n}{n!},\sq n\to\iy.$$
This gives finally $f(n)\sim 2$.
A complete asymptotic expansion follows by using more terms in the expansion of
$f(w)$ at $w=w_0$.

\subsect{Generalizing Kummer's identity}%
Kummer's identity reads:
$$
\F ab{1+a-b}{-1}=
\br{\G(1+a-b)\,\G\(1+\br12a\)}{\G(1+a)\,\G\(1+\br12a-b\)}.
$$
where $1+a-b\ne 0,-1,-2,\ldots$\ . 

A generalization (\ref{Vidunas} (2001)) is
written in the form
$$
\F{a+n}b{1+a-b}{-1}=
P(n)\br{\G(a-b)\,\G\(\br12a+\br12\)} {\G(a)\,\G\(\br12a+\br12-b\)}+
Q(n)\br{\G(a-b)\,\G\(\br12a\)} {\G(a)\,\G\(\br12a-b\)}.
$$
Vidunas showed that, for $n=-1, 0, 1, 2,\ldots, (a)_n\ne0$ and
$a-b\ne0,-1,-2,\ldots$:
$$\el{
P(n)&=    \br1{2^{n+1}}\FFF{-\br12n}{-\br12n-\br12}{\br12a-b}{\br12}{\br12a}1
=\br12\FFF{-\br12n}{-\br12n-\br12}{b}{-n}{\br a2}1,\cr
Q(n)&=\br{n+1}{2^{n+1}}
\FFF{-\br12n+\br12}{-\br12n}{\br12a+\br12-b}{\br32}{\br12a+\br12}1
=\br12\FFF{-\br12n+\br12}{-\br12n}{b}{-n}{\br a2+\br12}1.\cr
}\en{lf2}$$

We are interested in the asymptotic behaviour of $P(n)$ and $Q(n)$ for large
values of $n$.
Using, if $\Re d>\Re c> 1$,
$$
\br{(c)_k}{(d)_k}=
\br{\G(d)}{\G(d-c)\G(c)}\int_0^1t^{d-c-1}(1-t)^{c+k-1}\,dt,
$$
and a  few manipulations of the Gauss functions, we obtain the integral
representations, for $\Re a>\Re 2b>0$,
$$\el{
P(n)&=\br{2^{-n}\,\G(\br12a)}{\G(b)\,\G(\br12a-b)}
\intp\sinh^{a-2b-1}(t)\,\cosh^{-a-n}(t)\,\cosh(n+1)t\,dt,\cr
Q(n)&=\br{2^{-n}\,\G(\br12a+\br12)}{\G(b)\,\G(\br12a+\br12-b)}
\sq\sq\intp\sinh^{a-2b-1}(t)\,\cosh^{-a-n}(t)\,\sinh(n+1)t\,dt,\cr
}
\en{lf3}$$
We can use standard methods for obtaining the asymptotic behaviour of these
integrals. First we write the hyperbolic functions
$\cosh(n+1)t$ and
$\sinh(n+1)t$ as exponential functions.

We obtain for $n\to\iy$
$$\el{
P(n)&\sim\br{2^{2b-1}\G(\br12a)}{\G(b)\,\G(\br12a-b)}\sk
c_k\br{\G(b+\br12k)}{n^{b+\br12k}},\cr
Q(n)&\sim\br{2^{2b-1}\G(\br12a+\br12)}{\G(b)\,\G(\br12a+\br12-b)}\sk
c_k\br{\G(b+\br12k)}{n^{b+\br12k}},\cr
}$$
where 
$$c_0=1,\sq c_1=0,\sq c_2=\brg12(5b-4a+3).$$

These expansions are not  valid for $b=\br12a$, $b=\br12a+\br12$, respectively. 
In these cases we have
$$P(n)=\br1{2^{n+1}},\sq Q(n)=\br{n+1}{2^{n+1}},$$
as follows from \req{lf2}. Similar for the cases
$$b=\br12a+m, \sq b=\br12a+\br12+m,\sq m=0,1,2,\ldots.$$

For negative values of $n$ different integral representations should be used.
We have
$$\el{
P(-n-1)&=   2^n\br{\ph{1-\br12a}}{\ph{1-b}}
\FFF{-\br12n}{-\br12n+\br12}{\br12a-b}{\br12}{\br12a-n}1,\cr
Q(-n-1)&=   -n\,2^n\br{\ph{\br12-\br12a}}{\ph{1-b}}
\FFF{-\br12n+\br12}{-\br12n+1}{\br12a+\br12-b}{\br32}{\br12a+\br12-n}1.\cr
}$$
with the integral representations
$$\el{
P(-n-1) &= \ \  \br{ 2^{n+1}\,\G(1-b)}{\G(\br12a-b)\G(1-\br12a)} 
\int_0^{\br12\pi}\cos^{n+1-a} t\ \sin^{2c-1}t\ \cos nt\,dt,\cr
Q(-n-1) &=  -\br{ 2^{n+1}\,\G(1-b)}{\G(\br12a+\br12-b)\G(\br12-\br12a)} 
\int_0^{\br12\pi}\cos^{n+1-a} t\ \sin^{2c-1}t\ \sin nt\,dt.\cr
} $$
By integrating on contours $(i\infty, 0) \cup
(0,\br12\pi)\cup(\br12\pi,\br12\pi+i\infty)$, similar integrals arise as in
\req{lf3}, and again standard methods can be used for obtaining the asymptotic
behaviour.
\sect{ Concluding remarks}%
{\parindent=30pt
\item{[1.]\sq}
The asymptotic analysis of all 26 cases in
$$
\F{a+e_1\la}{b+e_2\la}{c+e_3\la}{z}, \sq e_j=0,\pm1,$$
can be reduced to the 4 cases
$$\vbox{\hsize=9.0truecm
\noindent
$$\vbox{\offinterlineskip
\halign{\strut             \vrule
\hfil $\ #\quad$\hfil   &  \vrule 
\hfil $\ #$\hfil   &       \vrule 
\hfil $\ #$ \hfil   &      \vrule
\hfil $\ #$ \hfil     \vrule 
    \cr \noalign{\hrule}
&\hfil \quad e_1 \quad\hfil  &  \hfil\quad e_2\quad\hfil  &  \hfil \quad
e_3\quad\hfil    \cr 
\noalign{\hrule}
\   &   &  \      &  \     \cr 
A & 0  &  0  &  +    \cr 
B & 0  &  -  &  +    \cr 
C & +  &  -  &  0    \cr 
D & +  &  2+  &  0   \cr 
\noalign{\hrule}
}} $$ 
}$$
\item{[2.]\sq}
These cases are of interest for orthogonal polynomials and special functions.
\item{[3.]\sq}
New recent results on uniform asymptotic expansions have been published or
announced for special cases.
\item{[4.]\sq}
The distribution of the zeros of Jacobi polynomials for non-classical values of
the parameters $\al$ and $\be$ shows interesting features.
New research is needed for describing the asymptotics of these distributions.
\item{[5.]\sq}
The asymptotics of $_3F_2$ terminating functions with large parameters
is quite difficult. Standard methods based on integrals and differential
equations are not available. Recursion relations may be explored further.
\par
}

\vfe
\sect{Bibliography}
\frenchspacing
\baselineskip=12pt
\parindent=30pt
\def\hfb{\hfill\break}
\def\ref#1{{\reffont#1\ }\advance \refnonum by 1}
\def\refno{{    {\the\refnonum}}}

\newcount\refnonum
\refnonum=1
\ninepoint
\item{[\refno]}\ref{M. Abramowitz and I.A. Stegun} (1964), 
{\sl Handbook of mathematical functions with formulas, graphs and mathematical
tables}, Nat. Bur. Standards Appl. Series, {\bf 55}, 
U.S. Government Printing Office, Washington, D.C. (paperback edition 
published by Dover, New York).

\item{[\refno]}
\ref{G.E. Andrews, R. Askey, R. Roy, Ranjan} (1999).
{\sl Special functions},
{Cambridge University Press}, Cambridge.

\item{[\refno]}
\ref{Christof Bosbach and Wolfgang Gawronski} (1999). 
Strong asymptotics for Jacobi polynomials with varying weights. 
{\sl Methods Appl. Anal.} 
{\bf 6},
39-54.

\item{[\refno]}
\ref{W.G.C. Boyd and T.M. Dunster} (1986). 
Uniform asymptotic solutions
of a class of second-order linear differential equations having a
turning point and a regular singularity, with an application to
Legendre functions, 
{\sl SIAM J. Math. Anal.}, 
{\bf 17}, 
422--450.

\item{[\refno]}
\ref{B.L. Braaksma} (1963).
{\sl Asymptotic expansions and analytic continuation 
for a class of Barnes-integrals}
Noordhoff, Groningen.

\item{[\refno]}
\ref{L.-C. Chen and M.E.H. Ismail} (1991). 
On asymptotics of Jacobi polynomials, 
{\sl SIAM J. Math. Anal.}, 
{\bf 22}, 1442--1449.

\item{[\refno]}
\ref{T.M. Dunster} (1991). 
Conical functions with one or both parameters large, 
{\sl Proc. Royal Soc. Edinburgh}, 
{\bf 119A}, 
311--327.

\item{[\refno]}
\ref{T.M. Dunster} (1999). 
Asymptotic approximations for the Jacobi 
and ultraspherical polynomials, and related functions, 
{\sl Methods and Applications of Analysis}, 
{\bf 6}, 
281--315.

\item{[\refno]}
\ref{T.M. Dunster} (2001).
Uniform asymptotic expansions for the reverse 
generalized Bessel polynomials, and related 
functions,
{\sl SIAM J. Math. Anal.} 
{\bf 32},
987--1013.

\item{[\refno]}
\ref{C.L. Frenzen and R. Wong} (1988).
Uniform asymptotic expansions of Laguerre polynomials,
{\sl SIAM J. Math. Anal.}, 
{\bf 19},  
1232--1248.

\item{[\refno]}
\ref{C.L. Frenzen} (1990). 
Error bounds for a uniform asymptotic expansion of the Legendre function 
$Q_n^{-m}(\cosh z)$,
{\sl SIAM J. Math. Anal.}, 
{\bf 21}, 
523--535.

\item{[\refno]}
\ref{Amparo Gil, Javier Segura and Nico M. Temme} 
(2000). 
Computing toroidal functions for wide ranges of the parameters,
{\sl J. Comput. Phys.}, 
{\bf 161}, 
204-217.

\item{[\refno]}
\ref{X.-S. Jin and R. Wong} (1998).
Uniform asymptotic expansions for Meixner polynomials,
{\sl Constructive Approximation}
{\bf 14}
113--150.

\item{[\refno]}
\ref{X.-S. Jin and R. Wong} (2001).
Asymptotic formulas for the zeros of the Meixner-Pollaczek polynomials
and their zeros,
{\sl Constr. Approx.}
{\bf 17}
59--90.

\item{[\refno]}
\ref{D.S. Jones} (2001).
Asymptotics of the hypergeometric function,
{\sl Math. Methods Appl. Sci.}
{\bf 24}
369--389.

\item{[\refno]}
\ref{R. Koekoek \& R.F. Swartouw}
\rm (1998).
The Askey-scheme of hypergeometric orthogonal 
polynomials and its $q-$analogue.
Technical University Delft. Report 98--17. 

\item{[\refno]}
\ref{X.-C. Li, and R. Wong} (2000).
A uniform asymptotic expansion for Krawtchouk polynomials,
{\sl J. Approximation Theory}
{\bf 106}
155--184.

\item{[\refno]}
\ref{X. Li and R. Wong} (2001).
On the asymptotics  of the Meixner polynomials,
{\sl J. Approximation Theory}
{\bf 96}
281--300.

\item{[\refno]}\ref{Y.L. Luke} (1969).
{\sl The special functions and their approximations},
Vol. I--II,
Academic Press, New York.

\item{[\refno]}
\ref{A. Mart\'\i nez-Finkelshtein, P. Mart\'\i nez-Gonz\'alez}  
\ref{and R. Orive} (2000).
Zeros of Jacobi polynomials with varying non-classical
              parameters
In {\sl Special functions} (Hong Kong, 1999), Charles Dunkl, Mourad Ismail
and Roderick Wong (eds.),
World Sci. Publishing, River Edge, NJ

\item{[\refno]}
\ref{A. Mart\'\i nez-Finkelshtein and N.M. Temme} (1999).
Private communications.

\item{[\refno]}
\ref{A. Olde Daalhuis} (2001).
Private communications.

\item{[\refno]}
\ref{A. Olde Daalhuis} (2001).
Uniform asymptotic expansions for hypergeometric 
functions with large parameters, I, II.
Submitted.

\item{[\refno]}
\ref{F.W.J. Olver}
\rm (1974 \&   1997).
\sl Asymptotics and Special Functions.
\rm Academic Press, New York. Reprinted in 1997 by A.K. Peters,

\item{[\refno]}
\ref{Bo Rui and R. Wong} (1994).
Uniform asymptotic expansion of Charlier polynomials,
{\sl Methods and Applications of Analysis}, 
{\bf 1}, 
294--313.

\item{[\refno]}
\ref{Bo Rui and R. Wong} (1996).
Asymptotic behavior of the Pollaczek polynomials and their zeros,
{\sl Stud. Appl. Math.}
{\bf 96}
307--338.

\item{[\refno]}
\ref{G. Szeg\"o} (1975).
{\sl Orthogonal polynomials}, 4th edition, Amer. Math.  Soc. 
Colloq. Publ. {\bf 23}, Providence, R.I.

\item{[\refno]}
\ref{N.M. Temme} (1986).
Uniform asymptotic expansion for a class of polynomials
bi\-orthogonal on the unit circle, 
{\sl Constr. Approx.}, 
{\bf 2}, 369--376.

\item{[\refno]}
\ref{N.M. Temme} (1990).
Asymptotic estimates for Laguerre polynomials,
{\sl ZAMP}, 
{\bf 41}, 
114--126.

\item{[\refno]}
\ref{N.M. Temme} (1996).
{\sl Special functions: An introduction to the classical functions
of mathematical physics},
Wiley, New York.

\item{[\refno]}
\ref{F. Ursell} (1984). 
Integrals with a large parameter: Legendre
functions of large degree and fixed order, {\sl Math. Proc. Camb.
Philos. Soc.}, 
{\bf 95}, 
367--380.

\item{[\refno]}
\ref{R. Vidunas} (2001).
A generalization of Kummer's identity.
Accepted for publication in {\sl Rocky Mountain J. Math.}.

\item{[\refno]}
\ref{G.N. Watson} (1918). 
Asymptotic expansions of hypergeometric functions, 
{\sl Trans. Cambridge Philos. Soc.}, 
{\bf 22}, 
277--308.

\item{[\refno]}
\ref{R. Wong} (1989).
{\sl Asymptotic approximations of integrals}, 
Academic Press, New York.

\item{[\refno]}
\ref{R. Wong and Q.-Q. Wang} (1992).
On the asymptotics of the Jacobi function and its zeros,
{\sl SIAM J. Math. Anal.}
{\bf 23}
1637--1649.

\item{[\refno]}
\ref{R. Wong and J.-M. Zhang} (1996).
A uniform asymptotic expansion for the Jacobi polynomials 
with explicit remainder,
{\sl Appl. Anal.}
{\bf 61}
17--29.

\item{[\refno]}
\ref{R. Wong and J.-M. Zhang} (1997).
The asymptotics of a second solution 
to the Jacobi differential equation,
{\sl Integral Transforms Spec. Funct.}
{\bf 5}
287--308.

\item{[\refno]}
\ref{R. Wong and J.-M. Zhang} (1997).
Asymptotic expansions of the generalized Bessel polynomials,
{\sl J. Comput. Appl. Math.}
{\bf 85}
87--112.

\bye